\let\ORIlabel\label
\let\ORIrefstepcounter\refstepcounter
\AddToHook{package/hyperref/before}{%
   \let\label\ORIlabel
   \let\refstepcounter\ORIrefstepcounter
}
\documentclass[onefignum,onetabnum]{siamart171218}

\usepackage{tikz}
\usepackage{pgfplots}
\usepgfplotslibrary{groupplots}
\usetikzlibrary{calc}
\pgfplotsset{compat=1.18}

\definecolor{ExpColor}{HTML}{0072B2}
\definecolor{PolyColor}{HTML}{E69F00}
\definecolor{LogColor}{HTML}{009E73}
\definecolor{UniformColor}{HTML}{D55E00}
\definecolor{OneClusterColor}{HTML}{CC79A7}
\definecolor{TwoClusterColor}{HTML}{4D4D4D}

\definecolor{c0}{HTML}{176B9B}
\definecolor{c1}{HTML}{BD5C25}
\definecolor{c2}{HTML}{397D65}
\definecolor{c3}{HTML}{8B5E9D}
\definecolor{c4}{HTML}{C24963}
\definecolor{c5}{HTML}{786936}

\usepackage{lipsum}
\usepackage{amsfonts}
\usepackage{graphicx}
\usepackage{epstopdf}
\usepackage{algorithmic}
\usepackage{booktabs}
\ifpdf
  \DeclareGraphicsExtensions{.eps,.pdf,.png,.jpg}
\else
  \DeclareGraphicsExtensions{.eps}
\fi

\newsiamremark{remark}{Remark}
\newsiamremark{hypothesis}{Hypothesis}
\crefname{hypothesis}{Hypothesis}{Hypotheses}
\newsiamthm{claim}{Claim}

\headers{FT-Cholesky}
        {Manavalan and Sch\"onlieb}

\title{Finite-time Cholesky factorization}

\author{
Rahul Manavalan\thanks{Center for Mathematical Sciences, Lund University
(\email{rahul.manavalan@matstat.lu.se})}
\and 
Carola Bibiane Sch\"onlieb\thanks{DAMTP, University of Cambridge}
}

\usepackage{amsopn}

\makeatletter
\newcommand*{\addFileDependency}[1]{
  \typeout{(#1)}%
  \@addtofilelist{#1}%
  \IfFileExists{#1}{}{\typeout{No file #1.}}
}
\makeatother

\newcommand*{\myexternaldocument}[1]{%
    \externaldocument{#1}%
    \addFileDependency{#1.tex}%
    \addFileDependency{#1.aux}%
}

\usepackage{tikz}
\usetikzlibrary{positioning}
\usepackage{subcaption}
\usepackage{booktabs}
\usepackage{multirow}
\usepackage{siunitx}

\ifpdf
\hypersetup{
  pdftitle={FT-Cholesky},
  pdfauthor={RMa}
}
\fi

\myexternaldocument{ex_supplement}

\begin{document}

\maketitle

\begin{abstract}
We propose an algorithm for estimating the Cholesky factorization of positive, semi-definite matrices without pivoting. Our method has the same data access pattern as the canonical Cholesky method, and therefore is amenable to familiar hardware optimization techniques. 
Our method interprets the factorization task as a variational problem and in turn solves it using an induced dynamical system. We study the rate of convergence of this dynamical system to its stationary state for both positive definite and semi-definite matrices. Numerical experiments validate that our method is competitive with pivoted Cholesky on several practical examples. 
\end{abstract}



\section{Introduction}
\label{sec:introduction}

The Cholesky factorization is one of the most studied transformations in numerical linear algebra. For a symmetric positive semi-definite matrix $A \in \mathbb{R}^{n \times n}$, it takes the form
\begin{align}
    A = LL^\top,
\end{align}
where $L \in \mathbb{R}^{n \times n}$ is lower triangular with $L_{ii} \geq 0$. This simple transformation has ubiquitous applications in different disciplines of the natural sciences.

The first documented use of this factorization is due to Benoit \cite{ref1}, who published the work of (Commandant) Cholesky in solving linear least-squares problems arising in geodetic networks. Originally developed to solve normal systems of equations, the method later morphed into a solver for positive definite systems
\begin{align}
    Ax = b,
\end{align}
for some $x,b \in \mathbb{R}^n$, through the works of Fox, Wilkinson \cite{10.1093/qjmam/1.1.149}, and Turing \cite{10.1093/qjmam/1.1.287}.
The factorization also provides a way to simplify other problems. As pointed out in \cite{10.1093/qjmam/1.1.149}, it can be used to transform generalized eigenvalue problems of the form
\begin{align}
    Ax = \lambda Bx,
\end{align}
where $A$ is symmetric, $B \in \mathbb{R}^{n \times n}$ is positive definite, and $\lambda \in \mathbb{R}$. Let $B = LL^\top$ be the Cholesky factorization of $B$, and use the change of variable $x := L^{-\top}y$. We then have
\begin{align}
    \begin{gathered}
        Ax = \lambda LL^\top x, \\
        L^{-1} A L^{-\top} y = \lambda y.
    \end{gathered}
\end{align}
In other words, we manage to turn a generalized eigenvalue problem into an eigenvalue problem corresponding to a symmetric matrix. While this is rarely used today as a direct method, sparse approximations of $L$ are used to precondition iterative methods like LOBPCG \cite{Kressner2023}. Lastly, the Cholesky factorization readily allows us to calculate the determinant of $A = LL^\top$ through
\begin{align}
    \det(A) = \prod_{i=1}^n L_{ii}^2.
\end{align}

The ubiquity of the Cholesky transformation rests on a combination of computational efficiency, numerical reliability, and efficient implementations. For positive definite matrices, it provides a unique factorization at substantially lower computational cost than an orthogonalization-based method such as Householder QR. For a matrix of size $n \times n$, the Cholesky decomposition requires approximately $n^3/3$ floating-point operations, compared with approximately $4n^3/3$ for Householder QR. For $n = 10^4$, these counts are approximately $3.3 \times 10^{11}$ and $1.3 \times 10^{12}$, respectively. At a sustained computing capacity of $2$ TFlops, these operation counts correspond to approximately $\frac{1}{6}$ seconds for Cholesky and four times longer for QR. Actual runtimes depend on the implementation and hardware, but the difference grows with matrix size and can be significant when the factorization is repeated within an outer-loop application.
The appeal of the Cholesky transformation is not efficiency alone. When $A \succ 0$, its factorization is unique \cite{doi:10.1137/1.9780898718027.ch10}, and its floating-point computation is backward stable \cite{doi:10.1137/1.9781611971811.ch8}. The semidefinite case ($A \succeq 0$) can admit infinitely many such factorizations. However, for a permutation matrix $P$, the pivoted matrix $P^\top A P$ has a unique Cholesky factorization \cite{doi:10.1137/1.9780898718027.ch10}, and the pivoted Cholesky transformation can also be shown to be backward stable \cite{10.1093/oso/9780198535645.003.0010}. These properties are supported by reliable LAPACK \cite{anderson1999lapack} implementations, which use blocked algorithms and BLAS-3 routines to compute the transformation efficiently; a brief summary is given in \cref{tab:cholesky-implementations}. Together, these three factors explain why the Cholesky method is so widely used.

\begin{table}[ht]
\centering
\caption{LAPACK Cholesky factorization routines, where
$x \in \{s,d,c,z\}$ denotes single-precision real, double-precision real,
single-precision complex, and double-precision complex arithmetic,
respectively. For a complete specification of these routines, see \cite{anderson1999lapack}.}
\label{tab:cholesky-implementations}
\begin{tabular}{ll}
\toprule
Routine & Description \\
\midrule
\texttt{xPOTRF}  & Dense positive-definite matrix; blocked algorithm \\
\texttt{xPOTF2}  & Dense positive-definite matrix; unblocked algorithm \\
\hline 
\texttt{xPSTRF}  & Positive-semidefinite matrix; complete pivoting \\
\texttt{xPSTF2}  & Positive-semidefinite matrix; unblocked complete pivoting \\
\bottomrule
\end{tabular}
\end{table}

\subsection{Previous work}
While pivoted Cholesky is an established method for factorizing positive semi-definite matrices; 
its implementation is more involved than the canonical Cholesky method. This can be attributed to the additional $\mathcal{O}(n^2)$ cost incurred for finding the next pivot; and the memory overhead for reordering the arrays during pivoting.
To avoid this, it is common practice to modify the positive-semi definite (PSD) matrix such that the canonical Cholesky algorithm can be employed while incurring some error due to the modification. 

One common modification is \textit{diagonal loading}. It entails shifting the spectrum of a PSD matrix $A$ by a factor $\mu \in \mathbb{R}_+$ such that $A + \mu I$ is positive definite. This inexpensive modification is popular in spatial statistics (a.k.a adding a nugget \cite{matheron1963principles}), in non-convex optimization (a.k.a Levenberg-Marquardt stabilization \cite{nocedal2006numerical}), in inverse problems (a.k.a Tikhonov regularization \cite{tikhonov1963regularization}) and in machine learning (a.k.a ridge regression \cite{scholkopf2002learning}). It is evident that the estimated Cholesky factor depends non-linearly on $\mu$. Therefore, finding an optimal shift parameter usually requires the specification of additional optimality criterion; whose choice is not unique. 

Another family of modifications is referred to as \textit{modified Cholesky methods} \cite{fang2008modified}. They originated in non-convex optimization, where Newton and quasi-Newton methods required a positive definite approximation of the Hessian at a given iteration. Gill, Murray and Wright \cite{gill1974newton} proposed strategies that entails finding a diagonal matrix $E \in \mathbb{R}^{n \times n}$ such that $A + E$ is positive definite, well-conditioned and sufficiently close to $A$. Schnabel and Eskow \cite{schnabel1990new} proposed another method that uses Gershgorin bounds. Cheng and Higham \cite{cheng1998modified} developed a different modification using a block $LDL^\top$ factorization. A systematic catalogue of modified Cholesky methods can be found in \cite{fang2008modified}. Unlike diagonal loading, modified Cholesky methods find litle use outside of optimization and numerical linear algebra. 

Orthogonal to these modification approaches, there exists another class of
methods that seek to eliminate the null space from a PSD matrix. These
so-called \textit{low-rank methods} approximate a matrix as
$A \approx L_r L_r^\top$, where $r \ll n$ is the rank of the approximation.
The simplest way to do this is to truncate the eigendecomposition of $A$.
This is more expensive than the Cholesky factorization itself, so it is seldom
used. On the other hand, Nystr\"om-based approaches \cite{drineas2005nystrom}, such as the randomly
pivoted Cholesky method \cite{chen2025randomly}, scale as $\mathcal{O}(r^2 n)$ and work extremely well
when $A$ has a low-rank structure. Furthermore, this family of randomized
methods is also stable for downstream applications involving approximations
of different matrix functions of these PSD matrices \cite{persson2023randomized, persson2025algorithmagnostic}.

Alternatively, one can mimic a pivoted Cholesky method through \textit{deflation} \cite{dostal2011cholesky}. Recall that a canonical Cholesky algorithm breaks down when the pivot, say $\alpha \leq 0$. In deflation, we ignore the update on this iteration; setting the corresponding column to zero. This method was primarily used to obtain an approximate Cholesky factor for the purpose of preconditioning linear systems arising from discretizations of PDEs. As noted in \cite{dostal2011cholesky}, the decision of when to deflate and when not to remains less understood, therefore deflation remains an approximate Cholesky factorization method. 

Lastly, Edelman and Holy \cite{edelman2026jordan} recently introduced a new method for factoring semi-definite and indefinite matrices using \textit{Jordan algebras} \cite{faraut1994analysis}. This method has the same structure of a canonical Cholesky algorithm except they use the hemiplex number system. A hemiplex number $z$ is of the form $z = a + m \mu + n \nu$ where $a,m,n \in \mathbb{R}$ and $1, \mu, \nu$ are the basis elements. The number system has additional axioms as specified in \cite{edelman2026jordan}.
In computing the Cholesky factorization of the following  matrix:
\begin{align}
    \begin{gathered}
        A = \begin{bmatrix}
            0 & 1 \\ 
            1 & 0
        \end{bmatrix}, \quad L = \begin{bmatrix}
            a & 0 \\ 
            b & c
        \end{bmatrix}; 
    \end{gathered}
\end{align}
one requires among other conditions that $a^2 = 0$ and $ab = 1$. This is impossible when $a,b \in \mathbb{R}$. On the other hand, the axioms of hemiplex numbers allow for this. Their implementations suggest that their method can be competitive with pivoted Cholesky and $LDL^\top$ decompositions. However, a comprehensive benchmark on practical examples remains wanting. 

Taken together, these approaches provide different ways to handle positive semidefinite matrices, either by modifying the matrix, constructing a low-rank approximation, handling vanishing pivots, or extending the underlying number system. This review has focused on approaches relevant to the factorization of dense PSD matrices, rather than providing an exhaustive account of sparse and scalable methods or their applications \cite{schafer2021sparse}. For a general introduction to Cholesky factorization of semidefinite matrices, its computation and numerical stability, we refer the reader to \cite{higham1990analysis} and references therein.

\subsection{Contributions}
\label{sec:contributions}
This work introduces a new method for calculating the Cholesky factorization of positive-semi definite matrices without pivoting. Its main contributions are as follows:
\begin{enumerate}
    \item We formulate Cholesky factorization as a variational problem; and derive dynamical systems whose stationary states provide the required optimality conditions. We derive analytical solutions to these dynamical systems and delineate how these solutions can be used to modify the canonical Choleksy algorithm and its blocked counterpart. We refer to the proposed modifications as \textit{finite-time Cholesky} methods. 
    \item We analyze the accuracy of the finite-time Cholesky map and discuss its dependence on the evolution time (of the dynamical system) $T$.
    \item We demonstrate with synthetic and practical examples that our methods is competitive with (the blocked and unblocked versions of) the pivoted Cholesky methods. 
\end{enumerate}

The remainder of the paper is organized as follows. The finite-time Cholesky algorithm is derived in \cref{sec:variational-problem} and \cref{sec:gradient-flow}. Its convergence rates with respect to the evolution time $T$ is studied in \cref{sec:dependence-on-T}.  Using synthetic matrices with different spectral shapes, the dependence of the accuracy of the method on $T$ is demonstrated in \cref{sec:numerical-experiments}. Lastly, we compare our method with pivoted Cholesky methods on different practical matrices in \cref{sec:applications}. 

\section{Finite time Cholesky factorization}
Let the space of positive-semi definite matrices be denoted by:
$$
\mathcal{A} := \{A \in \mathbb{R}^{n \times n} : A \succeq 0\}. 
$$
Similarly, the space of lower triangular matrices is denoted by:
$$
\mathcal{L}_n := \{L \in \mathbb{R}^{n \times n} : \forall i,j \in [1,n]^2, L_{ij} = 0, \textrm{ if } j >i\}.
$$
Recall that the latter is a linear space wheras the former is a cone in $\mathbb{R}^{n \times n}$. Any $A \in \mathcal{A}$ has a unique factorization $A =L L^\top$ where $L \in \mathcal{L}.$ This is referred to as the Cholesky factorization of $A$. For convenience, we will denote this using a map:
$$
L = \texttt{chol}(A),
$$
where $\texttt{chol} : \mathcal{A} \mapsto \mathcal{L}$ is a bijective map for all $A \succ 0$. On the other hand, when $A \succeq 0$, the factorization is not unique and one usually considers an equivalence class of Cholesky factors:
$$
[L]_A:= \{L \in \mathcal{L} : LL^T = A\}.  
$$

\subsection{Variational formulation}
\label{sec:variational-problem}

The Cholesky factorization of a given $A \in \mathcal{A}$ can be interpreted as a structured variational problem:
\begin{align}
    \begin{gathered}
       \arg \min_{L \in \mathcal{L}_n} \|A - LL^T\|_F^2. 
    \end{gathered}
    \label{eq:variational-problem}
\end{align}
Here, $\|\cdot\|_F$ is the Frobenius norm of a matrix. 
Denoting the matrices $A, L$ in block form such as:
\begin{align}
    \begin{gathered}
        A := \begin{bmatrix}
            \alpha & r^T \\ 
            r & B
        \end{bmatrix}, \quad 
        L := \begin{bmatrix}
            l & 0 \\ 
            u & \hat{L}
        \end{bmatrix},
    \end{gathered}
    \label{eq:block-matrix}
\end{align}
where $\alpha := A_{1,1}$, $r := A_{2:n, 1}$, $l := L_{1,1}$, $u := L_{2:n, 1}$ and $\hat{L} := L_{2:n, 2:n}$; 
one can show that the cost function in \cref{eq:variational-problem} can be decoupled. That is:

\begin{align}
    \begin{gathered}
        J[l,u,\hat{L}] := (l^2 - \alpha)^2 + 2 \|lu -r\|^2 + \|B - u u^T - \hat L \hat L^T \|_F^2. 
    \end{gathered}
\end{align}
This dependence of $J$ on $\hat L$ can be further dropped by examining its factorization in terms of its Schur complement.
\begin{align}
    \begin{gathered}
        A = \begin{bmatrix}
            1 & 0 \\ 
            \frac{r}{\alpha} & I 
        \end{bmatrix} \begin{bmatrix}
            \alpha & \\ 
            & B - \frac{1}{\alpha} r r^T
        \end{bmatrix} \begin{bmatrix}
            1 & \frac{r^T}{\alpha} \\ 
            0 & I 
        \end{bmatrix}. 
    \end{gathered}
\end{align}
We will denote one of the Schur complements by $S := B - \frac{1}{\alpha}r r^T$. Since, $A \succeq 0 \implies S \succeq 0$ and $uu^T = \frac{1}{\alpha} r r^T$, we can show that $\hat{L} \hat{L}^T = S$ which should be recursively computed by solving the variational problem $
        \arg \min_{\hat L \in \mathcal{L}_{n-1}} \|S - \hat L \hat L^T\|_F^2.$
Hence,
\begin{align}
    \begin{gathered}
        J[l,u,\hat L] = \hat J[l,u] =  (l^2 - \alpha)^2 + 2 \|lu -r\|^2. 
    \end{gathered}
\end{align}

The minimizer to this reduced functional can be obtained by solving the zeros of:
\begin{align}
    \begin{gathered}
        \partial_l \hat J = 4l(l^2-\alpha) + 4u^T(lu-r) = 0, \\ 
        \partial_u \hat J = 4l (lu-r) = 0. 
    \end{gathered}
\end{align}
When $\alpha > 0$ and using the constraint on $\mathcal{L}$ that $l \geq 0$, the minimizers are obtained by substitution. 
\begin{align}
    \begin{gathered}
        l_{>0} := \sqrt{\alpha}, \quad u_{>0} := \frac{r}{\sqrt{\alpha}}. 
    \end{gathered}
\end{align}

Substituting back, one arrives at the familiar recursive formulation of the Cholesky factorization. 
\begin{align}
    \begin{gathered}
        L = \textrm{chol}(A) = \begin{bmatrix}
            \sqrt{A_{1,1}} & 0 \\ 
            \frac{A_{2:n,1}}{\sqrt{A_{1,1}}} & \textrm{chol} \left(A_{2:n, 2:n} - \frac{1}{A_{1,1}} A_{2:n, 1} A_{2:n, 1}^\top\right) 
        \end{bmatrix}. 
    \end{gathered}
\end{align}

For positive definite matrices $A \succ 0$, this procedure for estimating $L$ is backwards stable. However, positive-semi-definite matrices require additional pivoting before $\texttt{chol}$ can be used reliably. 

\subsection{Cholesky factorization as a gradient flow}
\label{sec:gradient-flow}
Finite-time Cholesky factorization does not resort to pivoting. Instead, it solves the local minimization problem using its corresponding normalized gradient flow. 
\begin{align}
    \begin{gathered}
        \frac{d}{dt}l(t) = l(t)(\alpha - l(t)^2) + u(t)^T(r - l(t)u(t)), \quad l(0) = l_0, \\ 
        \frac{d}{dt} u(t) = l(t)(r-l(t)u(t)), \quad u(0) = 0 \in \mathbb{R}^{n-1}.  
    \end{gathered}
\end{align}
Here, $l_0 \geq 0$ is chosen appropriately (as we will specify later). Since $u(0) = 0$ and $\lim_{t \mapsto \infty} u(t) = \frac{1}{\sqrt{\alpha}} r$ for $\alpha > 0$ and $\lim_{t \mapsto \infty} u(t) = 0$ when $\alpha = 0$, it is sensible to constrain the dynamics of $u$ along the direction of the residual. Henceforth, we assume that, there exists $g :[0,T] \mapsto \mathbb{R}$ such that $g(0) = 0$ and $u(t) = r g(t)$. The corresponding reduced system of ODEs is given by:

\begin{align}
    \begin{gathered}
        \frac{d}{dt} l(t) = l(t) (\alpha - l(t)^2) + r^T r \; g(t) (1 - l(t) g(t)),  \\ 
        \frac{d}{dt} g(t) = l(t) (1 - l(t) g(t)). 
    \end{gathered}
\end{align}

To the best of our knowledge, there does not exist an analytical solution to this coupled ODE. 
As solving ODEs using numerical integrators for each iteration in a Cholesky decomposition function call is expensive, we will instead consider a surrogate ODE that has an analytical solution. 

\begin{align}
    \begin{gathered}
    \frac{d}{dt} l(t) = l(t) (\alpha - l(t)^2) ,
    \qquad l(0) = l_0 \geq \alpha^2 \\ 
    \frac{d}{dt} g(t) = l(t) (1 - l(t) g(t)), \quad g(0) = 0. 
    \end{gathered}
\end{align}
Observe that the equation are decoupled.
The first scalar equation is the logistic equation.
Its solution is given by:
\begin{equation}
    l(t) = \sqrt{\frac{\alpha l_0^2} {l_0^2 + (\alpha - l_0^2) e^{-2\alpha t}}}.
\end{equation}
Given $l(t)$, the second scalar equation is the prototypical non-homogeneous linear differential equation. Its solution is:
\begin{align}
    \begin{gathered}
    g(t) = \frac{1 - e^{-l(t)^2 t}}{l(t)}.
    \end{gathered}
\end{align}

Using these analytical solutions, one can define the finite-time Cholesky factorization of the matrix $A$ in \cref{eq:block-matrix} as follows:
\begin{equation}
    L(t) := \texttt{FT-chol}(A) = \begin{bmatrix}
        l(t) & 0 \\ 
        g(t) r & \texttt{FT-chol}(B - g(t)^2 r r^T)
    \end{bmatrix}.
\end{equation}

Notice that this surrogate ODE has the same stationary points as the reduced ODE. The only caveat of deleting the feedback term from its vector field is that $g(t)l(t) = 1$ no longer holds. We can argue that this effect is not really necessary. For one, when $A \succ 0$, $\lim_{t \mapsto \infty} g(t) l(t) = 1$. Furthermore when $A \succeq 0$, $r^Tr = 0$. Therefore, for computing the coefficients in the Cholesky decomposition, the surrogate ODE should suffice.

The issue of how the terminal time $t=T$ is chosen has to be taken up on a case-by-case basis. 
The initial condition $l_0$ on the other hand is simply chosen as $l_0 = \max\{1, \alpha^2\}$. 

The entire finite-time Cholesky factorization algorithm is specified in \cref{alg:finite-time-cholesky}. A blocked implementation of the same algorithm is listed in \cref{alg:blocked-finite-time-cholesky}. As the two implementations are equivalent in infinite precision arithmetic, it suffices to analyze the properties of the former.

\begin{algorithm}[H]
\caption{FT-Cholesky factorization}
\label{alg:finite-time-cholesky}
\begin{algorithmic}[1]
\REQUIRE Symmetric PSD matrix $A\in\mathbb{R}^{n\times n}$,
finite horizon $T>0$.
\ENSURE Lower-triangular finite-time factor $L_T$

\STATE Define
\[
\varphi_1(z)
=
\begin{cases}
\dfrac{e^z-1}{z}, & z\neq 0,\\[1ex]
1, & z=0.
\end{cases}
\]

\STATE $L\gets 0$

\FOR{$k=1,\ldots,n$}

    \STATE
    $\alpha_k
    \gets
    A_{kk}
    -
    L_{k,1:k-1}L_{k,1:k-1}^{\top}$

    \STATE
    $r_k
    \gets
    A_{k+1:n,k}
    -
    L_{k+1:n,1:k-1}L_{k,1:k-1}^{\top}$

    \STATE
    $s_k\gets\max\{1,\alpha_k\}$

    \STATE
    $p_k
    \gets
    \left[
    s_k^{-1}e^{-2\alpha_kT}
    +
    2T\,\varphi_1(-2\alpha_kT)
    \right]^{-1}$

    \STATE
    $\ell_k\gets\sqrt{p_k}$

    \STATE
    $g_k
    \gets
    T\,\ell_k\,\varphi_1(-p_kT)$

    \STATE
    $L_{kk}\gets\ell_k$

    \STATE
    $L_{k+1:n,k}\gets g_k r_k$

\ENDFOR

\STATE $L_T\gets L$

\RETURN $L_T$
\end{algorithmic}
\end{algorithm}

\Cref{alg:finite-time-cholesky} has a computational complexity of $\mathcal{O}(n^3/3)$. In comparison to the canonical Cholesky method, the method requires evaluation of two additional scalar functions, namely $x \mapsto e^x$ and $x \mapsto (e^x - 1)/ x$ per iteration. These are standard functions in the GNU scientific library \cite{galassi2009gsl} and can be included in existing Cholesky implementations without adding much computational overhead. \Cref{alg:blocked-finite-time-cholesky} has the same computational complexity. However, it reorganizes \cref{alg:finite-time-cholesky} to use BLAS-3 kernels \cite{dongarra1990blas3}. Therefore, it utilizes the underlying hardware more efficiently, particular when $n$ is large.

We will now demonstrate these algorithms on two positive-semi-definite matrix systems, namely the Hilbert matrices (denoted by $H \in \mathbb{R}^{n \times n}$) and the periodic finite-difference 1D-Laplacian matrices $D_{\Delta} \in \mathbb{R}^{n \times n}$. They are defined as follows:
\begin{align}
    \begin{gathered}
        [H]_{ij} = \frac{1}{i+j-1}, \quad \forall i,j \in \{1,2,\dots,n\}, \\ 
        [D_{\Delta}]_{ij} = \begin{cases}
            2, \quad i=j, \\ 
            -1, \quad j = i + 1 \mod n, \\ 
            0, \quad \textrm{else}. 
        \end{cases} \quad \forall i,j \in \{1,2,\dots,n\}. 
    \end{gathered}
\end{align}
The former is (strictly positive definite but) notoriously ill-conditioned for large $n$, and the latter satisfies $D_{\Delta} \mathbf{1} = 0$. Furthermore, in finite precision arithmetic, the representations of these matrices can also have negative eigenvalues. Promptly, one would expect canonical Cholesky routines like 
\texttt{xPOTF2} and \texttt{xPOTRF} to terminate prematurely. The finite time routines, on the other hand, return Cholesky factors that accurately reconstruct the matrices. This is shown in \cref{tab:examples}. 

\begin{algorithm}[H]
\caption{Block FT-Cholesky factorization}
\label{alg:blocked-finite-time-cholesky}
\begin{algorithmic}[1]
\REQUIRE Symmetric PSD matrix $A\in\mathbb{R}^{n\times n}$,
finite horizon $T>0$, block size $b\geq 1$.
\ENSURE Lower-triangular finite-time factor $L_T$

\STATE Define
\[
\varphi_1(z)
=
\begin{cases}
\dfrac{e^z-1}{z}, & z\neq 0,\\[1ex]
1, & z=0.
\end{cases}
\]

\STATE $L\gets 0$
\STATE $S\gets A$

\FOR{$k=1,1+b,\ldots,1+b\lfloor(n-1)/b\rfloor$}

    \STATE $q\gets\min(b,n-k+1)$
    \STATE $I\gets k:k+q-1$
    \STATE $J\gets k+q:n$

    \STATE $B\gets S_{I,I}$
    \STATE $C\gets S_{J,I}$

    \STATE $U\gets 0_{q\times q}$
    \STATE $V\gets 0_{|J|\times q}$

    \FOR{$j=1,\ldots,q$}

        \STATE
        $\alpha_j
        \gets
        B_{jj}
        -
        U_{j,1:j-1}U_{j,1:j-1}^{\top}$

        \STATE
        $r_j
        \gets
        B_{j+1:q,j}
        -
        U_{j+1:q,1:j-1}U_{j,1:j-1}^{\top}$

        \STATE
        $s_j
        \gets
        C_{:,j}
        -
        V_{:,1:j-1}U_{j,1:j-1}^{\top}$

        \STATE
        $\sigma_j\gets\max\{1,\alpha_j\}$

        \STATE
        $p_j
        \gets
        \left[
        \sigma_j^{-1}e^{-2\alpha_jT}
        +
        2T\,\varphi_1(-2\alpha_jT)
        \right]^{-1}$

        \STATE
        $\ell_j\gets\sqrt{p_j}$

        \STATE
        $g_j
        \gets
        T\,\ell_j\,\varphi_1(-p_jT)$

        \STATE $U_{jj}\gets\ell_j$
        \STATE $U_{j+1:q,j}\gets g_j r_j$
        \STATE $V_{:,j}\gets g_j s_j$

    \ENDFOR

    \STATE $L_{I,I}\gets U$
    \STATE $L_{J,I}\gets V$

    \STATE
    $S_{J,J}\gets S_{J,J}-VV^{\top}$

\ENDFOR

\STATE $L_T\gets L$

\RETURN $L_T$
\end{algorithmic}
\end{algorithm}

\begin{table}[H]
\centering
\caption{Relative Frobenius reconstruction errors for the Hilbert matrix
of order $14$ and the periodic one-dimensional finite-difference negative
Laplacian with $14$ grid points on $[0,1)$ and $h=1/14$, using Float64
arithmetic. FT-Cholesky uses $T=5\times10^{13}$ and block size $b=8$.
The reconstruction error is defined as $\|A-LL^\top\|_F/\|A\|_F$.}
\label{tab:examples}
\begin{tabular}{lrr}
\toprule
Method & Hilbert matrix & Periodic Laplacian \\
\midrule
FT-Cholesky, unblocked & $1.69\times10^{-14}$ & $1.10\times10^{-16}$ \\
FT-Cholesky, blocked   & $1.65\times10^{-14}$ & $1.10\times10^{-16}$ \\
\texttt{DPOTRF}        & --- & --- \\
\texttt{DPOTF2}        & --- & --- \\
\bottomrule
\end{tabular}
\end{table}

\subsection{Dependence on $T$}
\label{sec:dependence-on-T}

Note that in \cref{tab:examples}, the time evolution time $T$ of the finite-time algorithms was chosen arbitrarily as $T = 5 \times 10^{13}$. The dependence of the reconstruction accuracy on $T$ is not immediately obvious. This is studied next.  We will first consider positive definite matrices. 

\begin{theorem}[Positive definite matrices]
Let \(A\in \mathcal{A}\) with $A \succ 0$, and let \(L=\operatorname{chol}(A)\) denote its canonical Cholesky factor. For \(T>0\), let \(L_T=\texttt{FT-Chol}_T(A)\). Then,  \(L_T\) is well defined for every \(T>0\). Moreover, there exist constants \(C,C',\eta>0\) depending on \(A\), such that
\[
\|L_T-L\|_F\le Ce^{-\eta T},
\qquad
\|L_TL_T^\top-A\|_F\le C'e^{-\eta T},
\qquad T> 0.
\]
\end{theorem}

\begin{proof}
We will first establish arguments for the well-defined-ness of the Cholesky map. This is followed by a slight digression that bounds the error induced by the finite-time evolution. Lastly, we will establish the exponential rate of convergence of the finite time Cholesky map to the canonical Cholesky map using induction. These arguments are detailed below. \\ 
\textit{Step 1: Well-definedness.}  
Recall that the finite time Cholesky map is given by:
\begin{align*}
    \begin{gathered}
        L_T := \texttt{FT-chol}(A) = \begin{bmatrix}
            l(T) & 0 \\ 
            g(T) r & \texttt{FT-chol}(B - g(T)^2 r r^T)
        \end{bmatrix}. 
    \end{gathered}
\end{align*}  
To prove well-defined-ness, it suffices to show  that the Schur complement: $B - g(T)^2 r r^\top \succ 0$ for all $T > 0$. To prove this, consider the difference between the finite-time Schur complement and its canonical counterpart. 
\[
        S_T - S  = \left(\frac{1}{\alpha} - g(T)^2\right) rr^\top.
\]
Since, $S \succ 0$, it follows that if $\frac{1}{\alpha} > g(T)^2$, then $S_T \succ 0$. To show this, regroup the activation functions as:
\[
\frac{1}{l(T)^2} = \frac{l_0^2 + (\alpha - l_0^2) e^{-2\alpha T}}{\alpha l_0^2} = \frac{1}{\alpha} + \left(\frac{1}{l_0^2} - \frac{1}{\alpha}\right) e^{-2\alpha T}.
\]
Since, $l_0^2 \geq a$ by the construction of the FT-Cholesky algorithm, it follows that $1/l(T)^2 \leq 1/\alpha$. Recalling that for $\alpha > 0$: 
\[
g(T)^2 = \frac{1}{l(T)^2} (1- e^{- l(T)^2 T})^2, \quad 0 < 1- e^{- l(T)^2 T} < 1, \quad T>0; 
\]
we can conclude that $g(T)^2 < \frac{1}{\alpha}.$ Therefore, $L_T$ is well defined for all $T > 0$.  
\\ \\ 
\textit{Step 2: Bounds on scalar errors.} Bounds for the absolute difference between the scalar finite-time and canonical Cholesky coefficients are derived below. Note that:
\begin{align*}
\begin{gathered}
    \sqrt{\alpha} - l(T) = \frac{\alpha - l^2(T)}{\sqrt{\alpha} + l(T)} \\  = \sqrt{\alpha} \left(\frac{\alpha}{l_0^2} - 1\right) e^{-2 \alpha T} \frac{1}{\sqrt{1+ \left(\frac{\alpha}{l_0^2} - 1\right) e^{-2 \alpha T}}} \frac{1}{\sqrt{1+ \left(\frac{\alpha}{l_0^2} - 1\right) e^{-2 \alpha T}} + 1 }.
    \end{gathered}
\end{align*}
 Since, $l_0^2 \geq \alpha$; $\phi(\alpha, l_0) := \left| \frac{\alpha}{l_0^2} - 1\right| = (1 - \frac{\alpha}{l_0^2})$. It follows that:

\begin{align*}
    \begin{gathered}
        |\sqrt{\alpha} - l(T)| = \sqrt{\alpha} \; \phi(\alpha, l_0) e^{-2 \alpha T} \frac{1}{\sqrt{1 - \phi(\alpha, l_0) e^{-2 \alpha T}}} \frac{1}{\sqrt{1 - \phi(\alpha, l_0) e^{-2 \alpha T}}+1}, \\ 
        \leq \sqrt{\alpha} \; \phi(\alpha, l_0)\frac{1}{\sqrt{1 - \phi(\alpha, l_0)}} \frac{1}{\sqrt{1-\phi(\alpha, l_0)}+1}e^{-2 \alpha T}. 
    \end{gathered}
\end{align*}
Therefore,  $|\sqrt{\alpha} - l(T)| = \mathcal{O}(e^{-2 \alpha T})$. Similarly,
\begin{align*}
    \begin{gathered}
    \left| \frac{1}{\sqrt{\alpha}} - g(T) \right| \leq   \left| \frac{l(T) - \sqrt{\alpha}}{\sqrt{\alpha} \; l(T)} \right| +  \left| \frac{e^{- l(T)^2}}{l (T)} \right|. 
     \end{gathered}
\end{align*}
If $l_0^2 \geq \alpha$, then $l(T) \geq l_0 \geq \sqrt{\alpha}$; it follows that:
\begin{align*}
    \begin{gathered}
        \left| \frac{l(T) - \sqrt{\alpha}}{\sqrt{\alpha} \; l(T)} \right| \leq \frac{|l(T) - \sqrt{\alpha}|}{\sqrt{\alpha}} = \mathcal{O}(e^{-2 \alpha T}), \quad  
        \left| \frac{e^{- l(T)^2}}{l (T)} \right| \leq \frac{e^{- \alpha T}}{\sqrt{\alpha}} = \mathcal{O}(e^{- \alpha T }).
    \end{gathered}
\end{align*}
As a result, $|\frac{1}{\sqrt \alpha} - g(T)| = \mathcal{O}(e^{-\alpha T})$. 
\\ \\ 
\textit{Step 3: Convergence rate.} As specified in the proof sketch, we will derive convergence rates using induction. We will deduce the hypothesis for the convergence rates by examining two base cases, namely $n=2,3$.  
For a matrix given matrix $A_2 \succ 0 \in \mathbb{R}^{2 \times 2}$, we denote its canonical and finite time Cholesky factors by $[L]_2$ and $[L_T]_2$ respectively. They are given by:
\begin{align}
    \begin{gathered}
        A_2 := \begin{bmatrix}
            a_{11} & a_{12} \\ 
            a_{12} & a_{22}
        \end{bmatrix}, \quad [L]_2 := \begin{bmatrix}
            \sqrt {a_{11}} & 0 \\ 
            \frac{1}{\sqrt{a_{11}}} a_{12} & \sqrt{a_{22} - \frac{1}{a_{11}} a_{12}^2}
        \end{bmatrix}, \\ 
        [L_T]_2 :=  \begin{bmatrix}
            l_T(a_{11}) & 0 \\ 
            g_T(a_{11}) a_{12} & \sqrt{a_{22} - g_T^2(a_{11}) a_{12}^2}
        \end{bmatrix}. 
    \end{gathered}
    \label{eq:case-2}
\end{align}
Denoting the Schur complements by $s, s_T$, it follows that:
\begin{align*}
    \begin{gathered}
        \|[L]_2 - [L_T]_2\|_F^2 =  (\sqrt{a_{11}} - l_T(a_{11}))^2  + a_{12}^2 \left(\frac{1}{\sqrt{a_{11}}} - g_T(a_{11})\right)^2 +  (\sqrt{s} - \sqrt{s_T})^2 . 
    \end{gathered}
\end{align*}
We may bound the last term as follows:
\begin{align*}
    \begin{gathered}
        s - s_T = a_{12}^2 \left(g_T^2(a_{11}) - \frac{1}{a_{11}}\right) = a_{12}^2 \left(g_T(a_{11}) - \frac{1}{\sqrt{a_{11}}}\right) \left(g_T(a_{11}) + \frac{1}{\sqrt{a_{11}}}\right). 
    \end{gathered}
\end{align*}
Using $g_T(a_{11}) \leq \frac{1}{\sqrt{a_{11}}}$, $s_T \geq s$ and the previous scalar error bound, it follows that:
\begin{align*}
    \begin{gathered}
        |s - s_T| \leq 2a_{12}^2 \frac{1}{a_{11}} e^{- a_{11} T} \implies | \sqrt{s} - \sqrt{s_T}| = \frac{|s - s_T|}{\sqrt{s} + \sqrt{s_T}} \leq \frac{a_{12}^2}{a_{11} \sqrt{s}} e^{-a_{11} T}.   
    \end{gathered}
\end{align*}
Using this, we can show that \[\|[L]_2 - [L_T]_2\|_F = \mathcal{O}(e^{-a_{11}T}), \quad \|A_2 - [L_T]_2 [L_T]_2^\top\|_F = \mathcal{O}(e^{-a_{11}T}).\]
Similarly, let $A_3 \in \mathbb{R}^{3 \times 3}$ be symmetric
positive definite, and let $[L]_3$ and $[L_T]_3$ denote its
canonical and finite-time Cholesky factors, respectively.
Define the exact and finite-time Schur-complement entries by
\[
    s_{ij}
    := a_{ij} - \frac{a_{1i}a_{1j}}{a_{11}},
    \qquad
    s_{ij,T}
    := a_{ij} - g_T(a_{11})^2 a_{1i}a_{1j},
    \qquad i,j \in \{2,3\}.
\]
Then, the corresponding Cholesky factors are given by:
\begin{align*}
\begin{gathered}
A_3=\begin{bmatrix}
a_{11}&a_{12}&a_{13}\\
a_{12}&a_{22}&a_{23}\\
a_{13}&a_{23}&a_{33}
\end{bmatrix},\quad 
[L]_3=\begin{bmatrix}
\sqrt{a_{11}}&0&0\\
a_{12}/\sqrt{a_{11}}&\sqrt{s_{22}}&0\\
a_{13}/\sqrt{a_{11}}&s_{23}/\sqrt{s_{22}}&
\sqrt{s_{33}-s_{23}^2/s_{22}}
\end{bmatrix},\\[2pt]
[L_T]_3=\begin{bmatrix}
l_T(a_{11})&0&0\\
g_T(a_{11})a_{12}&l_T(s_{22,T})&0\\
g_T(a_{11})a_{13}&g_T(s_{22,T})s_{23,T}&
\sqrt{s_{33,T}-g_T(s_{22,T})^2s_{23,T}^2}
\end{bmatrix}.
\label{eq:case-3}
\end{gathered}
\end{align*}
Proceeding in a similar fashion as before, we can show that:
\begin{align*}
    \begin{gathered}
        \|[L_T]_3 - [L]_3\|_F = \mathcal{O}(e^{- \mu_3 T}), \\ \|A_3 - [L_T]_3[L_T]_3^{\top}\|_F = \mathcal{O}(e^{- \mu_3 T}), \quad \mu_3 = \min \left\{a_{11}, a_{22} - \frac{a_{12}^2}{a_{11}} \right\}. 
    \end{gathered}
\end{align*}
Note that: $l_{11}^2 = a_{11}$, $l_{22}^2 = a_{22} - \frac{a_{12}^2}{a_{11}}.$
This constitutes the base case for the induction argument. 
Assume that for analogously defined positive definite matrix $A_{k-1}$ and its Cholesky factors $[L]_{k-1}$ and $[L_T]_{k-1}$ for some $k-1 < n$, the following is true:
\begin{align}
\boxed{
    \begin{gathered}
        \|[L_T]_{k-1} - [L]_{k-1}\|_F = \mathcal{O}(e^{- \mu_{k-1} T}), \\
        \|A_{k-1} - [L_T]_{k-1}[L_T]_{k-1}^{\top}\|_F = \mathcal{O}(e^{- \mu_{k-1} T}), \\
        \mu_{k-1} = \min \left(l_{11}^2, \dots, l_{(k-2)(k-2)}^2\right).
    \end{gathered}
}
\end{align}
Denoting $A_k \in \mathbb{R}^{k \times k}$ and its finite time Cholesky factor $[L_T]_k$ by:
\begin{align*}
    \begin{gathered}
        A_k = \begin{bmatrix}
            a_{11} & a_{1, 2:k} \\
            a_{2:k,1} & a_{2:k, 2:k}
        \end{bmatrix}, \quad
        [L_T]_k = \begin{bmatrix}
            l_T(a_{11}) & 0 \\
            g_T(a_{11}) \; a_{2:k,1} & \texttt{FT-Chol}(a_{2:k, 2:k})
        \end{bmatrix},
    \end{gathered}
\end{align*}
the error in the factorization is given by:
\begin{align*}
    \begin{gathered}
        \|L_k - [L_T]_k \|_F^2 =
        (l_T(a_{11}) - \sqrt{a_{11}})^2
        + \|(g_T(a_{11}) - (1/\sqrt{a_{11}})) a_{2:k,1}\|^2
        + \\
        \|\texttt{Chol}(s_{2:k, 2:k}) -
        \texttt{FT-Chol}(s_{2:k, 2:k, T})\|_F^2.
    \end{gathered}
\end{align*}
Since the $\texttt{Chol} : \mathcal{A} \mapsto \mathcal{L}$ is Lipschitz continuous say with constant $K > 0$, it follows that:
\begin{align*}
    \begin{gathered}
        \|\texttt{Chol}(s_{2:k, 2:k}) -
        \texttt{FT-Chol}(s_{2:k, 2:k, T})\|_F^2
        \leq K \| s_{2:k, 2:k} - s_{2:k, 2:k, T} \|_F^2 \\
        + \|\texttt{Chol}(s_{2:k, 2:k, T}) -
        \texttt{FT-Chol}(s_{2:k, 2:k, T})\|_F^2.
    \end{gathered}
\end{align*}
We will first bound the first term to the right using the scalar bounds derived for $A_2$.
\begin{align*}
    \begin{gathered}
        s_{2:k, 2:k} - s_{2:k, 2:k, T}
        = \left(\frac{1}{a_{11}} - g_T(a_{11})^2 \right)
        a_{2:k,1} a_{2:k,1}^\top \\
        \implies
        \|s_{2:k, 2:k} - s_{2:k, 2:k, T}\|_F^2
        \leq \frac{2}{a_{11}} e^{- a_{11} T}
        \|a_{2:k, 1}\|^4_2.
    \end{gathered}
\end{align*}
From step 1, we know that $s_{2:k, 2:k, T} \succ 0, \forall T > 0$. Using the induction hypothesis, it follows that:
\begin{align*}
    \begin{gathered}
        \|\texttt{Chol}(s_{2:k, 2:k}) -
        \texttt{FT-Chol}(s_{2:k, 2:k, T})\|_F
        = \mathcal{O}(e^{-\nu_{k-1} T}).
    \end{gathered}
\end{align*}
Here, $\nu_{k-1}$ is the minimum of the square of the diagonal entries of $\texttt{Chol}(s_{2:k, 2:k, T})$. This quantity can in-turn be related to $l_{ii}$ where $i \in \{2,\dots,k\}$.
Firstly, $s_{2:k, 2:k, T} \succeq s_{2:k, 2:k}$. Therefore,
$\textrm{diag}(\texttt{Chol}(s_{2:k, 2:k, T}))_j
\geq \textrm{diag}(\texttt{Chol}(s_{2:k, 2:k}))_j$
for $j \in \{1,\dots, k-1\}$.
Additionally,
$
\textrm{diag}(\texttt{Chol}(s_{2:k, 2:k}))_j = l_{j+1,j+1}.
$
Consequently,
\begin{align*}
    e^{-\nu_{k-1} T} \leq e^{- \eta_{k-1} T},
    \quad
    \eta_{k-1} = \min \{l_{22}^2, \dots, l_{kk}^2\}.
\end{align*}
Therefore, for some constants $C_1, C_2 > 0$:
\begin{align*}
    \begin{gathered}
        \|L_k - [L_T]_k \|_F
        = C_1 e^{-l_{11}^2 T}
        + C_2 e^{-\eta_{k-1} T}
        = \mathcal{O}(e^{- \mu_k T}), \\
        \implies
        \|A_k - [L_T]_k [L_T]_k^\top\|_F
        = \mathcal{O}(e^{- \mu_k T}).
    \end{gathered}
\end{align*}
This completes the inductive argument. Hence the result is proven.
\end{proof}

\input{matern52_ft_cholesky_side_by_side}

This convergence rate is verified numerically for a Matern-5/2 matrix in \cref{fig:convergence-spd}.
We will now turn our attention to positive semi-definite matrices. Since, Cholesky factorization for such matrices over $\mathbb{R}$ is not unique, we will settle for the modest goal of bounding the growth of the reconstruction error with respect to $T$. 

\begin{theorem}[Positive semi-definite matrices]
    Let $A \in \mathcal{A}$. For $T > 0$, let $L_T = \texttt{FT-Chol}_T(A)$ denote its finite-time Cholesky factorization. Then, $L_T$ is well-defined for all $T > 0$. Moreover, there exists a constant $\bar C  >0$ such that for some $T_0 > 0$ that depends on $A$, such that:
    \[
    \|A - L_T L_T^\top\|_F \leq \bar C T^{-1}, \quad T \geq T_0. 
    \]
\end{theorem}

\begin{proof}
    The proof borrows several elements from and has the same structure as Theorem 2.1.
    Its arguments are detailed below. \\ 
\textit{Step 1: Well-definedness.} Recall that the finite-time Cholesky map is well-defined if:
\begin{align*}
    \begin{gathered}
        S_T - S  = \left(\frac{1}{\alpha} - g(T)^2\right) rr^\top \succeq 0.
    \end{gathered}
\end{align*}
Since, $A \succeq 0 \implies S \succeq 0$, it reduces to verifying $1/\alpha  \geq g(T)^2 $. Here, there are two cases. When $\alpha > 0$, we know from the proof of Theorem 2.1 that $S_T \succ 0$. On the other hand, when $\alpha = 0 \implies r = 0$ and $g(T) < \infty$.
Therefore, $S_T = B$ which is positive semi-definite by definition. \\ 

\textit{Step 2: Bounds on scalar errors.} The activation functions $l, g$ have different asymptotic behavior as $\alpha \mapsto 0$. We can show that:
\begin{align*}
    \begin{gathered}
        \lim_{\alpha \mapsto 0} l(T) = \frac{l_0}{\sqrt{1 + 2 l_0^2 T }}, \quad \lim_{\alpha \mapsto 0} g(T) = \frac{\sqrt{1 + 2 l_0^2 T}}{l_0} \left(1 - \exp \left[- \frac{T l_0^2}{1 + 2l_0^2 T}\right] \right). 
    \end{gathered}
\end{align*}
Using this, we can bound the scalar error functions appearing in Step 3. Firstly,  $\forall \alpha > 0$:
\begin{align*}
    \begin{gathered}
        \alpha - l_T(\alpha)^2  =  \frac{(\alpha - l_0^2) e^{-2 \alpha T}}{\frac{l_0^2}{\alpha} + (1 - \frac{l_0^2}{\alpha}) e^{-2 \alpha T}}. 
    \end{gathered}
\end{align*}
Since, $l_0^2 \geq \alpha$ by construction; the denominator of this expression  $\geq 1$, therefore:
\begin{align*}
    \begin{gathered}
        |\alpha - l_T^2(\alpha)| \leq |\alpha - l_0^2| e^{-2 \alpha T} = \mathcal{O}(e^{-2 \alpha T}). 
    \end{gathered}
\end{align*}
We are also interested in the case $\alpha = 0$. 
\begin{align*}
    \begin{gathered}
        \left| \lim_{\alpha \mapsto 0} \alpha - l_T^2(\alpha) \right| = \frac{l_0^2}{1+2T l_0^2} = \mathcal{O}(T^{-1}). 
    \end{gathered}
\end{align*}
Secondly, it will be of interest to obtain the bounds for $e^{- T l_T(\alpha)^2)}$. Since, $l_T^2(\alpha) \leq \alpha$, it follows that $e^{- l_T^2(\alpha) T } \leq e^{- \alpha T}.$ In other words, $e^{-  l_T(\alpha)^2 T)} = \mathcal{O}(e^{- \alpha T)}.$ \\ 

\textit{Step 3: Convergence rate.} The reconstruction error rate of $A$ will be obtained by induction. We use the definition in $\cref{eq:case-2}$. Firstly, we observe that for $A_2 \in \mathbb{R}^{2 \times 2} \succeq 0$, 
\begin{align*}
    \begin{gathered}
        \|A_2 - [L_T]_2 [L_T]_2^\top\|_F^2 = (a_{11} - l_T(a_{11})^2 )^2 + 2 a_{12}^2 e^{-2 l_T^2(a_{11}) T} + (s_T - l_T^2(s_T))^2. 
    \end{gathered}
\end{align*}
Furthermore when $s_T > 0$,
\begin{align}
    \begin{gathered}
        |s_T - l_T^2(s_T)| =  \frac{l_0^2 - s_T}{1 + l_0^2 \frac{e^{2 s_T T} - 1}{s_T}}. 
    \end{gathered}
\end{align}
Here, we slightly modify the algorithm and maintain a running update of $l_0$. 
In this case $l_0$ is $\max \{1, s_{T}, a_{11}\}$. 
Since $e^{2 s_T T} - 1 > 2s_T T$, it follows that:
\begin{align*}
    \begin{gathered}
       |s_T - l_T^2(s_T)| \leq \frac{l_0^2 - s_T}{1 + 2 Tl_0^2} \leq \frac{l_0^2}{1 + 2 Tl_0^2} = \mathcal{O}(T^{-1}). 
    \end{gathered}
\end{align*}
Similarly, when $s_T$ approaches zero:
\begin{align*}
    \begin{gathered}
        \left| \lim_{s_T \mapsto 0} s_T - l_T^2(s_T) \right| = \frac{l_0^2}{1 + 2 T l_0^2} = \mathcal{O}(T^{-1}).  
    \end{gathered}
\end{align*}
Assuming that the modified $l_0$ is used in both evaluations of the activation function $l$,
it follows that:
\begin{align*}
    \begin{gathered}
        \|A_2 - [L_T]_2 [L_T]_2^\top\|_F \leq  \frac{2 l_0^2}{1 + 2T l_0^2} + \sqrt2 |a_{12}|  e^{- a_{11} T}   = \mathcal{O}(T^{-1}). 
    \end{gathered}
\end{align*}
Similarly, for $A_3 \ \succeq 0$ and its corresponding Cholesky factor $[L_T]_3$, we can use the definitions in \cref{eq:case-3}, proceed in an analogous fashion and show that:
\begin{align*}
    \begin{gathered}
        \|A_3 - [L_T]_3 [L_T]_3^\top\|_F \leq \frac{3 l_0^2}{1 + 2 T l_0^2} + \sqrt{2(a_{12}^2 + a_{13}^2)} \; e^{- a_{11} T} + \sqrt 2 |s_{23, T}| e^{-s_{22,T} T}. 
    \end{gathered}
\end{align*}
Here, $l_0 = \max \{1, a_{11}, s_{22,T}, s_{33,T} - g_T^2(s_{22, T}) s_{23,T}^2\}$. 
The last term is dependent on $T$ and  can be bounded using the following:
\begin{align*}
    \begin{gathered}
        e^{-s_{22,T} T} \leq e^{-s_{22} T }, \quad  
        s_{23,T} \geq s_{23}, \quad
        |s_{23}| \leq \sqrt{s_{33} s_{22}} \\ 
        s_{22} \leq a_{22}, \quad s_{33} \leq a_{33}, \quad s_{22} \geq 0, \\ 
        \implies |s_{23, T}| e^{-s_{22,T} T} \leq \sqrt{a_{22} a_{33}} e^{-s_{22} T } \leq \sqrt{a_{22} a_{33}}, \quad T \geq 0.  
    \end{gathered}
\end{align*}
Therefore, $ \|A_3 - [L_T]_3 [L_T]_3^\top\|_F = \mathcal{O}(T^{-1})$. 
These two examples constitute the base cases for the induction hypothesis. For $A_k \succeq 0$, we hypothesize that: 
\begin{align}
\boxed{
\begin{gathered}
\left\|A_k-[L_T]_k[L_T]_k^\top\right\|_F
 = \mathcal{O}(T^{-1}).
\end{gathered}
}
\end{align}
Then, the proposition can be proven for $A_{k+1} \succeq 0$ as follows. We first define:
\begin{align*}
    \begin{gathered}
        A_{k+1} := \begin{bmatrix}
            A_{k} & a \\ 
            a^\top & a_{k+1, k+1}
        \end{bmatrix}, \quad 
        [L_T]_{k+1} = \begin{bmatrix}
            [L_T]_{k} & 0 \\ 
            \lambda_T^\top  & (L_T)_{k+1, k+1}
        \end{bmatrix}. 
    \end{gathered}
\end{align*}
Here, $\lambda_T = (L_T)_{k+1,1:k}$ and $a \in \mathbb{R}^{k}$. Then,
\begin{align*}
    \begin{gathered}
        \|A_{k+1} - [L_T]_{k+1} [L_T]_{k+1}^\top\|_F^2 = \|A_k-[L_T]_k[L_T]_k^\top\|_F^2 + 2 \|a - [L_T]_k \lambda_T\|^2  \\ + \;  (a_{k+1, k+1} - \lambda_T^\top \lambda_T - (L_T)_{k+1, k+1}^2 )^2.   
    \end{gathered}
\end{align*}
Denoting by $r_T := a_{k+1, k+1} - \lambda_T^\top \lambda_T$ and choosing $l_0 \geq a_{k+1, k+1}$ note that:
\begin{align*}
    \begin{gathered}
        (a_{k+1, k+1} - \lambda_T^\top \lambda_T - (L_T)_{k+1, k+1}^2 )^2 = |r_T - l_T^2(r_T)|^2 \leq \left( \frac{l_0^2}{1 + 2 Tl_0^2}\right)^2. 
    \end{gathered}
\end{align*}
Lastly, it suffices  to show that $\|a - [L_T]_k \lambda_T\|^2$ is bounded. 
\begin{align*}
    \begin{gathered}
        \|a - [L_T]_k \lambda_T\|^2 \leq \|a\|^2 + 2 \|a\| \|[L_T]_k \lambda_T\| + \|[L_T]_k \lambda_T\|^2
    \end{gathered}
\end{align*}
To this end, we use the fact that $r_T \geq 0$ to show that
\begin{align*}
    \begin{gathered}
        \|\lambda_T\|^2 = \sum_{i=1}^k (L_T)_{k+1,i}^2 \leq a_{k+1, k+1}.
\end{gathered}
\end{align*}
Furthermore, defining by $x_i := a_{ii} - \sum_{j=1}^{i-1} ([L_T]_k)_{ij}^2$ we can show that
\begin{align*}
\begin{gathered}
        \|[L_T]_k\|_F^2 = \textrm{tr}([L_T]_k [L_T]_k^\top) = \sum_{i=1}^k ([L_T]_k [L_T]_k^\top)_{ii} =  \sum_{i=1}^k \sum_{j=1}^i ([L_T]_k)_{ih}^2, \\ 
        \sum_{j=1}^i ([L_T]_k)_{ih}^2 = \sum_{j=1}^{i-1} ([L_T]_k)_{ih}^2 + ([L_T]_k)_{ii}^2 \leq a_{ii} + (x_i - l_T(x_i)) \leq a_{ii} + \frac{l_0^2}{1 + 2 T l_0^2}, \\ 
        \implies \|[L_T]_k\|_F^2 \leq \textrm{tr}(A_k) + \frac{k l_0^2}{1 + 2 T l_0^2}. 
    \end{gathered}
\end{align*}
Taken together, we have shown that $\|a - [L_T]_k \lambda_T\|^2 = \mathcal{O}(T^{-1})$. Using this and the induction hypothesis $\|A_{k+1} - [L_T]_{k+1} [L_T]_{k+1}^\top\| = \mathcal{O}(T^{-1})$. 
\end{proof}

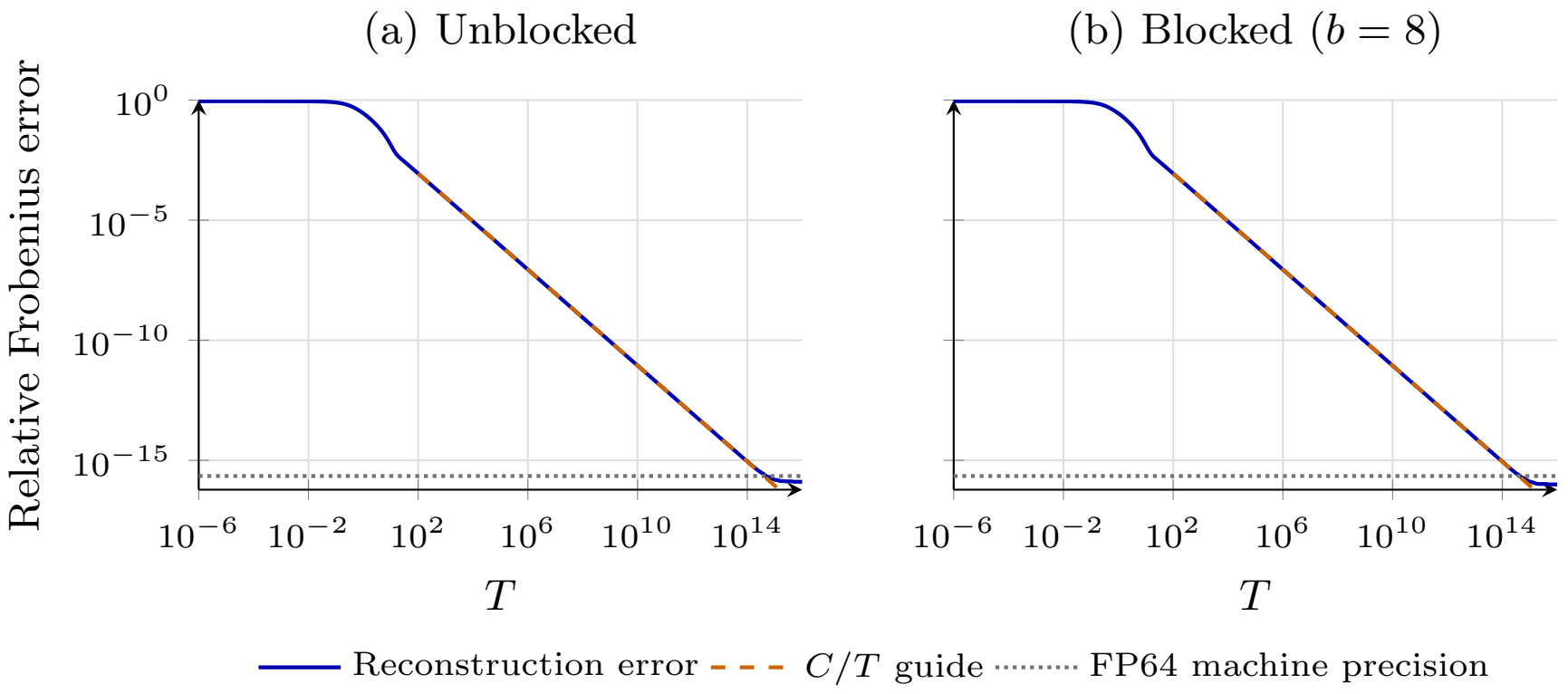
\begin{figure}[htbp]
\setlength{\abovecaptionskip}{2pt}
\centering
\pgfplotstableread{
T unblocked blocked
1.0000000000e-06 8.954292388931e-01 8.954292388931e-01
1.2589254118e-06 8.954290070426e-01 8.954290070426e-01
1.5848931925e-06 8.954287151599e-01 8.954287151599e-01
1.9952623150e-06 8.954283477010e-01 8.954283477010e-01
2.5118864315e-06 8.954278850969e-01 8.954278850969e-01
3.1622776602e-06 8.954273027118e-01 8.954273027118e-01
3.9810717055e-06 8.954265695307e-01 8.954265695307e-01
5.0118723363e-06 8.954256465077e-01 8.954256465077e-01
6.3095734448e-06 8.954244844863e-01 8.954244844863e-01
7.9432823472e-06 8.954230215814e-01 8.954230215814e-01
1.0000000000e-05 8.954211798824e-01 8.954211798824e-01
1.2589254118e-05 8.954188613039e-01 8.954188613039e-01
1.5848931925e-05 8.954159423597e-01 8.954159423597e-01
1.9952623150e-05 8.954122675842e-01 8.954122675842e-01
2.5118864315e-05 8.954076412484e-01 8.954076412484e-01
3.1622776602e-05 8.954018169300e-01 8.954018169300e-01
3.9810717055e-05 8.953944843783e-01 8.953944843783e-01
5.0118723363e-05 8.953852529744e-01 8.953852529744e-01
6.3095734448e-05 8.953736309004e-01 8.953736309004e-01
7.9432823472e-05 8.953589989025e-01 8.953589989025e-01
1.0000000000e-04 8.953405772408e-01 8.953405772408e-01
1.2589254118e-04 8.953173840511e-01 8.953173840511e-01
1.5848931925e-04 8.952881828745e-01 8.952881828745e-01
1.9952623150e-04 8.952514165237e-01 8.952514165237e-01
2.5118864315e-04 8.952051237007e-01 8.952051237007e-01
3.1622776602e-04 8.951468338280e-01 8.951468338280e-01
3.9810717055e-04 8.950734343371e-01 8.950734343371e-01
5.0118723363e-04 8.949810031029e-01 8.949810031029e-01
6.3095734448e-04 8.948645967115e-01 8.948645967115e-01
7.9432823472e-04 8.947179826762e-01 8.947179826762e-01
1.0000000000e-03 8.945333003768e-01 8.945333003768e-01
1.2589254118e-03 8.943006311485e-01 8.943006311485e-01
1.5848931925e-03 8.940074522431e-01 8.940074522431e-01
1.9952623150e-03 8.936379418460e-01 8.936379418460e-01
2.5118864315e-03 8.931720922980e-01 8.931720922980e-01
3.1622776602e-03 8.925845752096e-01 8.925845752096e-01
3.9810717055e-03 8.918432839605e-01 8.918432839605e-01
5.0118723363e-03 8.909074543497e-01 8.909074543497e-01
6.3095734448e-03 8.897252304356e-01 8.897252304356e-01
7.9432823472e-03 8.882304966382e-01 8.882304966382e-01
1.0000000000e-02 8.863387350101e-01 8.863387350101e-01
1.2589254118e-02 8.839415840909e-01 8.839415840909e-01
1.5848931925e-02 8.808996705301e-01 8.808996705301e-01
1.9952623150e-02 8.770331609992e-01 8.770331609992e-01
2.5118864315e-02 8.721093625145e-01 8.721093625145e-01
3.1622776602e-02 8.658266511859e-01 8.658266511859e-01
3.9810717055e-02 8.577941965413e-01 8.577941965413e-01
5.0118723363e-02 8.475077291526e-01 8.475077291526e-01
6.3095734448e-02 8.343236745022e-01 8.343236745022e-01
7.9432823472e-02 8.174385213190e-01 8.174385213190e-01
1.0000000000e-01 7.958887369390e-01 7.958887369390e-01
1.2589254118e-01 7.685990561087e-01 7.685990561087e-01
1.5848931925e-01 7.345181487334e-01 7.345181487334e-01
1.9952623150e-01 6.928730024662e-01 6.928730024662e-01
2.5118864315e-01 6.435192647698e-01 6.435192647698e-01
3.1622776602e-01 5.872605170879e-01 5.872605170879e-01
3.9810717055e-01 5.259316208133e-01 5.259316208133e-01
5.0118723363e-01 4.621175696114e-01 4.621175696114e-01
6.3095734448e-01 3.986069597340e-01 3.986069597340e-01
7.9432823472e-01 3.378490114781e-01 3.378490114781e-01
1.0000000000e+00 2.816333594746e-01 2.816333594746e-01
1.2589254118e+00 2.310320996736e-01 2.310320996736e-01
1.5848931925e+00 1.865094104719e-01 1.865094104719e-01
1.9952623150e+00 1.480787374485e-01 1.480787374485e-01
2.5118864315e+00 1.154387793709e-01 1.154387793709e-01
3.1622776602e+00 8.809256837149e-02 8.809256837149e-02
3.9810717055e+00 6.548321575817e-02 6.548321575817e-02
5.0118723363e+00 4.711740235476e-02 4.711740235476e-02
6.3095734448e+00 3.259266270389e-02 3.259266270389e-02
7.9432823472e+00 2.155000823289e-02 2.155000823289e-02
1.0000000000e+01 1.364656424365e-02 1.364656424365e-02
1.2589254118e+01 8.529499563986e-03 8.529499563986e-03
1.5848931925e+01 5.692287236768e-03 5.692287236768e-03
1.9952623150e+01 4.271177222445e-03 4.271177222445e-03
2.5118864315e+01 3.407068027760e-03 3.407068027760e-03
3.1622776602e+01 2.733017805868e-03 2.733017805868e-03
3.9810717055e+01 2.182078553735e-03 2.182078553735e-03
5.0118723363e+01 1.738257120198e-03 1.738257120198e-03
6.3095734448e+01 1.383594110914e-03 1.383594110914e-03
7.9432823472e+01 1.100809113533e-03 1.100809113533e-03
1.0000000000e+02 8.755301713223e-04 8.755301713223e-04
1.2589254118e+02 6.961706870231e-04 6.961706870231e-04
1.5848931925e+02 5.534383230567e-04 5.534383230567e-04
1.9952623150e+02 4.398962154564e-04 4.398962154564e-04
2.5118864315e+02 3.496017187184e-04 3.496017187184e-04
3.1622776602e+02 2.778120255576e-04 2.778120255576e-04
3.9810717055e+02 2.207456080263e-04 2.207456080263e-04
5.0118723363e+02 1.753897175962e-04 1.753897175962e-04
6.3095734448e+02 1.393455678307e-04 1.393455678307e-04
7.9432823472e+02 1.107041476175e-04 1.107041476175e-04
1.0000000000e+03 8.794680882308e-05 8.794680882307e-05
1.2589254118e+03 6.986581452422e-05 6.986581452420e-05
1.5848931925e+03 5.550092093616e-05 5.550092093615e-05
1.9952623150e+03 4.408880835535e-05 4.408880835535e-05
2.5118864315e+03 3.502278993548e-05 3.502278993547e-05
3.1622776602e+03 2.782072965028e-05 2.782072965027e-05
3.9810717055e+03 2.209950962461e-05 2.209950962460e-05
5.0118723363e+03 1.755471787098e-05 1.755471787097e-05
6.3095734448e+03 1.394449414847e-05 1.394449414847e-05
7.9432823472e+03 1.107668593891e-05 1.107668593890e-05
1.0000000000e+04 8.798638290831e-06 8.798638290831e-06
1.2589254118e+04 6.989078690762e-06 6.989078690758e-06
1.5848931925e+04 5.551667886036e-06 5.551667886032e-06
1.9952623150e+04 4.409875164292e-06 4.409875164288e-06
2.5118864315e+04 3.502906408128e-06 3.502906408120e-06
3.1622776602e+04 2.782468854700e-06 2.782468854696e-06
3.9810717055e+04 2.210200760892e-06 2.210200760888e-06
5.0118723363e+04 1.755629403739e-06 1.755629403735e-06
6.3095734448e+04 1.394548866465e-06 1.394548866461e-06
7.9432823472e+04 1.107731344745e-06 1.107731344730e-06
1.0000000000e+05 8.799034227529e-07 8.799034227529e-07
1.2589254118e+05 6.989328512744e-07 6.989328512666e-07
1.5848931925e+05 5.551825514592e-07 5.551825514553e-07
1.9952623150e+05 4.409974621902e-07 4.409974621862e-07
2.5118864315e+05 3.502969162141e-07 3.502969162102e-07
3.1622776602e+05 2.782508449970e-07 2.782508449813e-07
3.9810717055e+05 2.210225744016e-07 2.210225743977e-07
5.0118723363e+05 1.755645167017e-07 1.755645166978e-07
6.3095734448e+05 1.394558812378e-07 1.394558812300e-07
7.9432823472e+05 1.107737620222e-07 1.107737620144e-07
1.0000000000e+06 8.799073824048e-08 8.799073823657e-08
1.2589254118e+06 6.989353497666e-08 6.989353497275e-08
1.5848931925e+06 5.551841279706e-08 5.551841279315e-08
1.9952623150e+06 4.409984567072e-08 4.409984567072e-08
2.5118864315e+06 3.502975438489e-08 3.502975436926e-08
3.1622776602e+06 2.782512410571e-08 2.782512409008e-08
3.9810717055e+06 2.210228242911e-08 2.210228242520e-08
5.0118723363e+06 1.755646744114e-08 1.755646742551e-08
6.3095734448e+06 1.394559807130e-08 1.394559806348e-08
7.9432823472e+06 1.107738248052e-08 1.107738247661e-08
1.0000000000e+07 8.799077790238e-09 8.799077774607e-09
1.2589254118e+07 6.989355994607e-09 6.989355990699e-09
1.5848931925e+07 5.551842870089e-09 5.551842866182e-09
1.9952623150e+07 4.409985570147e-09 4.409985554517e-09
2.5118864315e+07 3.502976058621e-09 3.502976058621e-09
3.1622776602e+07 2.782512804064e-09 2.782512796249e-09
3.9810717055e+07 2.210228492995e-09 2.210228485180e-09
5.0118723363e+07 1.755646914875e-09 1.755646910968e-09
6.3095734448e+07 1.394559919277e-09 1.394559915370e-09
7.9432823472e+07 1.107738314481e-09 1.107738310573e-09
1.0000000000e+08 8.799078286500e-10 8.799078130197e-10
1.2589254118e+08 6.989356385364e-10 6.989356346288e-10
1.5848931925e+08 5.551843026392e-10 5.551842948241e-10
1.9952623150e+08 4.409985714728e-10 4.409985558424e-10
2.5118864315e+08 3.502976320429e-10 3.502976281353e-10
3.1622776602e+08 2.782512870493e-10 2.782512714190e-10
3.9810717055e+08 2.210228481273e-10 2.210228481273e-10
5.0118723363e+08 1.755646914876e-10 1.755646836724e-10
6.3095734448e+08 1.394559958353e-10 1.394559802050e-10
7.9432823472e+08 1.107738431709e-10 1.107738275406e-10
1.0000000000e+09 8.799078950797e-11 8.799078560035e-11
1.2589254118e+09 6.989356776133e-11 6.989355213098e-11
1.5848931925e+09 5.551844472210e-11 5.551844081446e-11
1.9952623150e+09 4.409985597519e-11 4.409984815996e-11
2.5118864315e+09 3.502977492725e-11 3.502977101957e-11
3.1622776602e+09 2.782514316325e-11 2.782513925555e-11
3.9810717055e+09 2.210230317870e-11 2.210229927097e-11
5.0118723363e+09 1.755647579210e-11 1.755646016160e-11
6.3095734448e+09 1.394560231941e-11 1.394558668887e-11
7.9432823472e+09 1.107739096069e-11 1.107738705282e-11
1.0000000000e+10 8.799078951709e-12 8.799075043754e-12
1.2589254118e+10 6.989363420161e-12 6.989359512108e-12
1.5848931925e+10 5.551844473656e-12 5.551840565480e-12
1.9952623150e+10 4.409984427068e-12 4.409976611193e-12
2.5118864315e+10 3.502985310168e-12 3.502981401640e-12
3.1622776602e+10 2.782518226789e-12 2.782502595285e-12
3.9810717055e+10 2.210242044230e-12 2.210238135145e-12
5.0118723363e+10 1.755646411514e-12 1.755638594546e-12
6.3095734448e+10 1.394574695729e-12 1.394570785761e-12
7.9432823472e+10 1.107731288173e-12 1.107731285238e-12
1.0000000000e+11 8.799079042958e-13 8.799039929297e-13
1.2589254118e+11 6.989480762305e-13 6.989441638830e-13
1.5848931925e+11 5.551883694032e-13 5.551844558201e-13
1.9952623150e+11 4.410090113672e-13 4.410050962287e-13
2.5118864315e+11 3.503141842392e-13 3.503102671428e-13
3.1622776602e+11 2.782584944120e-13 2.782545748503e-13
3.9810717055e+11 2.210242407842e-13 2.210164111997e-13
5.0118723363e+11 1.755674221861e-13 1.755634956138e-13
6.3095734448e+11 1.394614347206e-13 1.394575032303e-13
7.9432823472e+11 1.107798441730e-13 1.107759064915e-13
1.0000000000e+12 8.799869682176e-14 8.799475134690e-14
1.2589254118e+12 6.990664520426e-14 6.990268992176e-14
1.5848931925e+12 5.552679668986e-14 5.552282905391e-14
1.9952623150e+12 4.411671343540e-14 4.411273026776e-14
2.5118864315e+12 3.502383273777e-14 3.502373990635e-14
3.1622776602e+12 2.783395333083e-14 2.781820302185e-14
3.9810717055e+12 2.211724512882e-14 2.211318679042e-14
5.0118723363e+12 1.755719996886e-14 1.754920323650e-14
6.3095734448e+12 1.395062659295e-14 1.394648003309e-14
7.9432823472e+12 1.107089522604e-14 1.107060154149e-14
1.0000000000e+13 8.804690076877e-15 8.788681055232e-15
1.2589254118e+13 6.995719896029e-15 6.991335924371e-15
1.5848931925e+13 5.561939934485e-15 5.545710983909e-15
1.9952623150e+13 4.417396879934e-15 4.412735162422e-15
2.5118864315e+13 3.519123481789e-15 3.514269526786e-15
3.1622776602e+13 2.789013629345e-15 2.772194124966e-15
3.9810717055e+13 2.230956403670e-15 2.225557085194e-15
5.0118723363e+13 1.763015253064e-15 1.757220886055e-15
6.3095734448e+13 1.412480763349e-15 1.406219333280e-15
7.9432823472e+13 1.113160544807e-15 1.102507295816e-15
1.0000000000e+14 9.077214712832e-16 9.001550723744e-16
1.2589254118e+14 7.148011885996e-16 7.062496293080e-16
1.5848931925e+14 5.731268819690e-16 5.518740605358e-16
1.9952623150e+14 4.595730275127e-16 4.449565450068e-16
2.5118864315e+14 3.740326821572e-16 3.612432609789e-16
3.1622776602e+14 3.124178417992e-16 2.867859593976e-16
3.9810717055e+14 2.668131971150e-16 2.504030297432e-16
5.0118723363e+14 2.108750186216e-16 1.948477812151e-16
6.3095734448e+14 2.016208980362e-16 1.810363156590e-16
7.9432823472e+14 1.730970336387e-16 1.496366501672e-16
1.0000000000e+15 1.559262593761e-16 1.287977131148e-16
1.2589254118e+15 1.494255973623e-16 1.227270798239e-16
1.5848931925e+15 1.371718088736e-16 1.109661638313e-16
1.9952623150e+15 1.371718088736e-16 1.038584360269e-16
2.5118864315e+15 1.332188605364e-16 1.038584360269e-16
3.1622776602e+15 1.332188605364e-16 1.038584360269e-16
3.9810717055e+15 1.332188605364e-16 1.038584360269e-16
5.0118723363e+15 1.279571861963e-16 9.935001693534e-17
6.3095734448e+15 1.279571861963e-16 9.935001693534e-17
7.9432823472e+15 1.279571861963e-16 9.935001693534e-17
1.0000000000e+16 1.279571861963e-16 9.935001693534e-17
}\maternNullityData
\begin{tikzpicture}
\begin{groupplot}[
  group style={group size=2 by 1,horizontal sep=0.085\linewidth},
  width=0.46\linewidth,height=0.34\linewidth,
  xmin=1e-6,xmax=1e16,
  ymin=6e-17,ymax=1,
  xmode=log,ymode=log,log basis x=10,log basis y=10,
  xtick={1e-6,1e-2,1e2,1e6,1e10,1e14},
  ytick={1e-15,1e-10,1e-5,1},
  minor tick num=0,
  grid=major,grid style={black!12},axis lines=left,
  tick label style={font=\scriptsize},
  title style={font=\small},
  label style={font=\small},
  legend style={font=\scriptsize,draw=none,fill=none},
  legend cell align=left,
  legend columns=3,
  clip=true,
]
\nextgroupplot[
  title={(a) Unblocked},
  xlabel={$T$},
  ylabel={Relative Frobenius error},
  legend to name=matern-nullity-legend,
]
\addplot[thick,blue!70!black]
  table[x=T,y=unblocked] {\maternNullityData};
\addlegendentry{Reconstruction error}
\addplot[thick,dashed,orange!80!black]
  coordinates {(1e2,8.799078222749e-04) (1e16,8.799078222749e-18)};
\addlegendentry{$C/T$ guide}
\addplot[thick,densely dotted,black!55]
  coordinates {(1e-6,2.220446049250313e-16)
               (1e16,2.220446049250313e-16)};
\addlegendentry{FP64 machine precision}

\nextgroupplot[
  title={(b) Blocked ($b=8$)},
  xlabel={$T$},
  yticklabel=\empty,
]
\addplot[thick,blue!70!black]
  table[x=T,y=blocked] {\maternNullityData};
\addplot[thick,dashed,orange!80!black]
  coordinates {(1e2,8.799078222749e-04) (1e16,8.799078222749e-18)};
\addplot[thick,densely dotted,black!55]
  coordinates {(1e-6,2.220446049250313e-16)
               (1e16,2.220446049250313e-16)};
\end{groupplot}
\path (group c1r1.south west) -- (group c2r1.south east)
  node[midway,below=26pt] {\pgfplotslegendfromname{matern-nullity-legend}};
\end{tikzpicture}
\caption{Relative reconstruction error
$\| A-L_TL_T^\top\|_F/\| A\|_F$ for a $64\times64$
Mat\'ern-$5/2$ covariance matrix of rank $54$. We use the first $54$
sites of an $8\times8$ grid on $[0,1]^2$, length scale $\ell=0.2$,
and a nugget of $10^{-8}$ before duplicating ten row--column pairs.
Both methods follow the $C/T$ guide, with
$C=\sqrt{10}/(2\| A\|_F)$, before reaching the Float64 rounding floor.}
\label{fig:psd-rate}
\end{figure}

This rate is verified for a Matern-covariance matrix with controlled null-space in \cref{fig:psd-rate}. Hence, we have shown through Theorem 2.1 and Theorem 2.2 that for positive semi-definite matrices the accuracy of the finite-time Cholesky algorithm increases monotonously in the evolution parameter $T$ so long as it is well-defined. 

%

\section{Numerical experiments}
\label{sec:numerical-experiments}

We will now use the finite-time Cholesky algorithm and its blocked implementation to factor a collection of synthetic matrices. These matrices are randomly generated and are designed to have different conditioning and spectral profiles. They are constructed as follows. 
Let $W \in \mathbb{R}^{n \times n}$ be a random matrix whose elements $W_{ij} \sim \mathcal{N}(0,1)$. The column space corresponding to $W$ is orthogonalized using a QR decomposition, such that $W = Q  R$. The synthetic matrix is defined as $
        A_s = Q \Lambda_s Q^\top.$  
Here, $\Lambda_s \in \mathbb{R}^{n \times n}$ is a diagonal matrix. Its smallest and largest non-zero elements say $[\Lambda_s]_{\min}, [\Lambda_s]_{\max}$ specify the condition number $\kappa$. 
Furthermore, $[\Lambda_s]_{ii} \in [0, \infty)$ where $i \in [1,n]$. The distribution of these elements on $n$ points on this interval, in turn, specifies the spectral shape of $A_s$. 
The spectral distributions considered in this work are specified in \cref{tab:spectral-shapes}. 

\begin{table}[t]
\centering
\caption{Spectral constructions
$\Lambda_s=\operatorname{diag}(\lambda_1^{(s)},\ldots,\lambda_n^{(s)})$.
Here $a=\kappa^{-1}$, $p=\log(\kappa)/\log(n)$, and $m=n/2$,
with even $n\geq4$.
The values $u_i$ are sorted independent $\mathcal{U}(0,1)$ samples,
with endpoints reset to $u_1=0$ and $u_n=1$.
}
\label{tab:spectral-shapes}
\begin{tabular}{lll}
\toprule
Spectrum $s$ & Eigenvalue $\lambda_i^{(s)}$ & Index range \\
\midrule
Exponential
& $a^{\,1-(i-1)/(n-1)}$
& $i=1,\ldots,n$ \\[1ex]
Polynomial
& $(i/n)^p$
& $i=1,\ldots,n$ \\[1ex]
Log-uniform
& $a^{\,1-u_i}$
& $i=1,\ldots,n$ \\[1ex]
Uniform
& $a+(1-a)u_i$
& $i=1,\ldots,n$ \\[1ex]
One-cluster
& $a$
& $i=1$ \\[1ex]
& $0.9+0.1(i-2)/(n-2)$
& $i=2,\ldots,n$ \\[1ex]
Two-cluster
& $a\bigl(1+0.1(i-1)/(m-1)\bigr)$
& $i=1,\ldots,m$ \\[1ex]
& $0.9+0.1(i-m-1)/(m-1)$
& $i=m+1,\ldots,n$ \\
\bottomrule
\end{tabular}
\end{table}

Our experiment entails estimating the finite-time Cholesky factorization, say $L_T$, of a given matrix $A$ and measuring the Frobenius norm relative residual :
\begin{align}
    \begin{gathered}
        \frac{\|A - L_T L_T^\top\|_F}{\|A\|_F}. 
    \end{gathered}
\end{align}

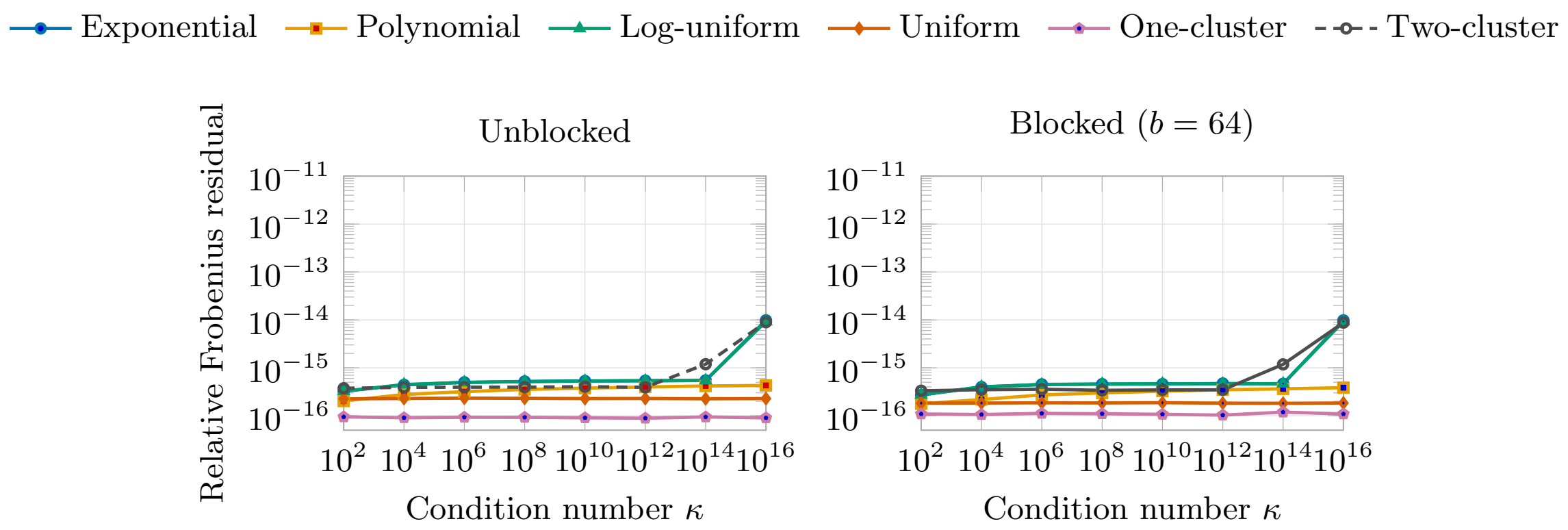
\begin{figure}[H]
\centering
\begin{tikzpicture}
\begin{groupplot}[
    group style={group size=2 by 1, horizontal sep=1.5cm},
    width=0.435\textwidth,
    height=0.31\textwidth,
    xmode=log,
    ymode=log,
    xmin=1e2, xmax=1e16,
    ymin=5e-17, ymax=1e-11,
    xtick={1e2,1e4,1e6,1e8,1e10,1e12,1e14,1e16},
    xticklabels={$10^2$,$10^4$,$10^6$,$10^8$,$10^{10}$,$10^{12}$,$10^{14}$,$10^{16}$},
    ytick={1e-16,1e-15,1e-14,1e-13,1e-12,1e-11},
    yticklabels={$10^{-16}$,$10^{-15}$,$10^{-14}$,$10^{-13}$,$10^{-12}$,$10^{-11}$},
    grid=major,
    major grid style={gray!25, line width=0.25pt},
    axis line style={gray!65},
    tick style={gray!65},
    tick label style={font=\small},
    label style={font=\small},
    title style={font=\small},
    xlabel={Condition number $\kappa$},
    scaled ticks=false,
    clip=true,
]
\nextgroupplot[
    title={Unblocked},
    ylabel={Relative Frobenius residual},
    legend to name=legendSPD0,
    legend columns=6,
    legend style={draw=none, fill=none, font=\small, /tikz/every even column/.append style={column sep=0.18cm}},
]
\addplot+[color=ExpColor, mark=*, mark size=1.25pt, line width=0.9pt] coordinates {(100, 3.24000000e-16) (10000, 4.40000000e-16) (1000000, 4.97000000e-16) (100000000, 5.19000000e-16) (10000000000, 5.31000000e-16) (1000000000000, 5.38000000e-16) (100000000000000, 5.46000000e-16) (10000000000000000, 9.79000000e-15)};
\addlegendentry{Exponential}
\addplot+[color=PolyColor, mark=square*, mark size=1.25pt, line width=0.9pt] coordinates {(100, 2.03000000e-16) (10000, 2.78000000e-16) (1000000, 3.18000000e-16) (100000000, 3.54000000e-16) (10000000000, 3.76000000e-16) (1000000000000, 3.98000000e-16) (100000000000000, 4.18000000e-16) (10000000000000000, 4.29000000e-16)};
\addlegendentry{Polynomial}
\addplot+[color=LogColor, mark=triangle*, mark size=1.25pt, line width=0.9pt] coordinates {(100, 3.26000000e-16) (10000, 4.48000000e-16) (1000000, 4.96000000e-16) (100000000, 5.15000000e-16) (10000000000, 5.30000000e-16) (1000000000000, 5.37000000e-16) (100000000000000, 5.45000000e-16) (10000000000000000, 9.28000000e-15)};
\addlegendentry{Log-uniform}
\addplot+[color=UniformColor, mark=diamond*, mark size=1.25pt, line width=0.9pt] coordinates {(100, 2.25000000e-16) (10000, 2.29000000e-16) (1000000, 2.33000000e-16) (100000000, 2.31000000e-16) (10000000000, 2.28000000e-16) (1000000000000, 2.29000000e-16) (100000000000000, 2.26000000e-16) (10000000000000000, 2.28000000e-16)};
\addlegendentry{Uniform}
\addplot+[color=OneClusterColor, mark=pentagon*, mark size=1.25pt, line width=0.9pt] coordinates {(100, 9.45000000e-17) (10000, 9.08000000e-17) (1000000, 9.28000000e-17) (100000000, 9.25000000e-17) (10000000000, 9.05000000e-17) (1000000000000, 8.89000000e-17) (100000000000000, 9.43000000e-17) (10000000000000000, 9.02000000e-17)};
\addlegendentry{One-cluster}
\addplot+[color=TwoClusterColor, mark=o, mark size=1.25pt, line width=0.9pt] coordinates {(100, 3.74000000e-16) (10000, 3.91000000e-16) (1000000, 3.94000000e-16) (100000000, 3.95000000e-16) (10000000000, 4.04000000e-16) (1000000000000, 3.92000000e-16) (100000000000000, 1.19000000e-15) (10000000000000000, 8.96000000e-15)};
\addlegendentry{Two-cluster}
\nextgroupplot[
    title={Blocked ($b=64$)},
    ylabel={},
]
\addplot+[color=ExpColor, mark=*, mark size=1.25pt, line width=0.9pt, forget plot] coordinates {(100, 2.71000000e-16) (10000, 3.97000000e-16) (1000000, 4.51000000e-16) (100000000, 4.56000000e-16) (10000000000, 4.57000000e-16) (1000000000000, 4.71000000e-16) (100000000000000, 4.63000000e-16) (10000000000000000, 9.79000000e-15)};
\addplot+[color=PolyColor, mark=square*, mark size=1.25pt, line width=0.9pt, forget plot] coordinates {(100, 1.78000000e-16) (10000, 2.18000000e-16) (1000000, 2.72000000e-16) (100000000, 2.99000000e-16) (10000000000, 3.26000000e-16) (1000000000000, 3.48000000e-16) (100000000000000, 3.64000000e-16) (10000000000000000, 3.85000000e-16)};
\addplot+[color=LogColor, mark=triangle*, mark size=1.25pt, line width=0.9pt, forget plot] coordinates {(100, 2.73000000e-16) (10000, 4.02000000e-16) (1000000, 4.46000000e-16) (100000000, 4.60000000e-16) (10000000000, 4.66000000e-16) (1000000000000, 4.63000000e-16) (100000000000000, 4.67000000e-16) (10000000000000000, 9.28000000e-15)};
\addplot+[color=UniformColor, mark=diamond*, mark size=1.25pt, line width=0.9pt, forget plot] coordinates {(100, 1.83000000e-16) (10000, 1.84000000e-16) (1000000, 1.86000000e-16) (100000000, 1.85000000e-16) (10000000000, 1.87000000e-16) (1000000000000, 1.83000000e-16) (100000000000000, 1.82000000e-16) (10000000000000000, 1.84000000e-16)};
\addplot+[color=OneClusterColor, mark=pentagon*, mark size=1.25pt, line width=0.9pt, forget plot] coordinates {(100, 1.09000000e-16) (10000, 1.06000000e-16) (1000000, 1.12000000e-16) (100000000, 1.10000000e-16) (10000000000, 1.07000000e-16) (1000000000000, 1.03000000e-16) (100000000000000, 1.19000000e-16) (10000000000000000, 1.09000000e-16)};
\addplot+[color=TwoClusterColor, mark=o, mark size=1.25pt, line width=0.9pt, forget plot] coordinates {(100, 3.32000000e-16) (10000, 3.49000000e-16) (1000000, 3.57000000e-16) (100000000, 3.40000000e-16) (10000000000, 3.46000000e-16) (1000000000000, 3.47000000e-16) (100000000000000, 1.18000000e-15) (10000000000000000, 8.78000000e-15)};
\end{groupplot}
\node[anchor=south] at ([yshift=2mm]current bounding box.north) {\pgfplotslegendfromname{legendSPD0}};
\end{tikzpicture}
\caption{FT-Cholesky reconstruction residuals for positive definite matrices. Left: unblocked FT-Cholesky; right: blocked FT-Cholesky with block size $b=64$. The relative Frobenius residual is plotted against the target condition number. Each point represents one realization.}
\label{fig:ftchol-trend-spd}
\end{figure}

We will consider one set of matrices from each of the groups in \cref{tab:spectral-shapes}. Each of the matrices in a group are specified the condition numbers $\kappa \in \{10^2, 10^4, 10^6, \dots, 10^{16}\}.$ The results of this experiments for one realization of these synthetic matrices is furnished in \cref{fig:ftchol-trend-spd}. 

Subsequently, we will modify the construction of these synthetic matrices to make them semidefinite. This is realized by randomly assigning certain eigenvalues in the spectrum to zero. 
We refer to the number of zero eigenvalues of a given matrix as its nullity. 
The reconstruction errors for nullities $1,2,5,10$ for the same synthetic matrices are specified in \cref{fig:ftchol-trend-psd1}, \cref{fig:ftchol-trend-psd2}, \cref{fig:ftchol-trend-psd5} and \cref{fig:ftchol-trend-psd10} respectively. 

\begin{figure}[H]
\centering
\begin{tikzpicture}
\begin{groupplot}[
    group style={group size=2 by 1, horizontal sep=1.15cm},
    width=0.435\textwidth,
    height=0.31\textwidth,
    xmode=log,
    ymode=log,
    xmin=1e2, xmax=1e16,
    ymin=5e-17, ymax=1e-11,
    xtick={1e2,1e4,1e6,1e8,1e10,1e12,1e14,1e16},
    xticklabels={$10^2$,$10^4$,$10^6$,$10^8$,$10^{10}$,$10^{12}$,$10^{14}$,$10^{16}$},
    ytick={1e-16,1e-15,1e-14,1e-13,1e-12,1e-11},
    yticklabels={$10^{-16}$,$10^{-15}$,$10^{-14}$,$10^{-13}$,$10^{-12}$,$10^{-11}$},
    grid=major,
    major grid style={gray!25, line width=0.25pt},
    axis line style={gray!65},
    tick style={gray!65},
    tick label style={font=\small},
    label style={font=\small},
    title style={font=\small},
    xlabel={Condition number $\kappa$},
    scaled ticks=false,
    clip=true,
]
\nextgroupplot[
    title={Unblocked},
    ylabel={Relative Frobenius residual},
    legend to name=legendPSD1,
    legend columns=6,
    legend style={draw=none, fill=none, font=\small, /tikz/every even column/.append style={column sep=0.18cm}},
]
\addplot+[color=ExpColor, mark=*, mark size=1.25pt, line width=0.9pt] coordinates {(100, 5.07000000e-16) (10000, 1.28000000e-14) (1000000, 1.23000000e-15) (100000000, 1.51000000e-15) (10000000000, 1.60000000e-15) (1000000000000, 2.80000000e-15) (100000000000000, 6.33000000e-15) (10000000000000000, 1.03000000e-14)};
\addlegendentry{Exponential}
\addplot+[color=PolyColor, mark=square*, mark size=1.25pt, line width=0.9pt] coordinates {(100, 2.96000000e-16) (10000, 5.31000000e-16) (1000000, 5.36000000e-16) (100000000, 1.51000000e-15) (10000000000, 3.76000000e-16) (1000000000000, 3.13000000e-15) (100000000000000, 7.45000000e-16) (10000000000000000, 8.78000000e-16)};
\addlegendentry{Polynomial}
\addplot+[color=LogColor, mark=triangle*, mark size=1.25pt, line width=0.9pt] coordinates {(100, 3.32000000e-16) (10000, 1.05000000e-15) (1000000, 5.03000000e-16) (100000000, 1.51000000e-15) (10000000000, 1.22000000e-15) (1000000000000, 2.08000000e-14) (100000000000000, 1.63000000e-15) (10000000000000000, 9.68000000e-15)};
\addlegendentry{Log-uniform}
\addplot+[color=UniformColor, mark=diamond*, mark size=1.25pt, line width=0.9pt] coordinates {(100, 7.15000000e-16) (10000, 3.19000000e-16) (1000000, 7.00000000e-16) (100000000, 6.97000000e-16) (10000000000, 2.33000000e-16) (1000000000000, 2.28000000e-16) (100000000000000, 4.05000000e-16) (10000000000000000, 8.49000000e-16)};
\addlegendentry{Uniform}
\addplot+[color=OneClusterColor, mark=pentagon*, mark size=1.25pt, line width=0.9pt] coordinates {(100, 3.74000000e-15) (10000, 3.68000000e-16) (1000000, 9.96000000e-17) (100000000, 2.99000000e-16) (10000000000, 1.34000000e-16) (1000000000000, 3.12000000e-14) (100000000000000, 1.52000000e-13) (10000000000000000, 5.86000000e-14)};
\addlegendentry{One-cluster}
\addplot+[color=TwoClusterColor, mark=o, mark size=1.25pt, line width=0.9pt] coordinates {(100, 1.39000000e-15) (10000, 5.16000000e-16) (1000000, 3.95000000e-16) (100000000, 3.95000000e-16) (10000000000, 4.03000000e-16) (1000000000000, 4.07000000e-16) (100000000000000, 1.23000000e-15) (10000000000000000, 8.83000000e-15)};
\addlegendentry{Two-cluster}
\nextgroupplot[
    title={Blocked ($b=64$)},
    ylabel={},
]
\addplot+[color=ExpColor, mark=*, mark size=1.25pt, line width=0.9pt, forget plot] coordinates {(100, 1.00000000e-15) (10000, 1.25000000e-14) (1000000, 1.25000000e-15) (100000000, 1.35000000e-15) (10000000000, 1.64000000e-15) (1000000000000, 4.71000000e-16) (100000000000000, 5.95000000e-15) (10000000000000000, 1.03000000e-14)};
\addplot+[color=PolyColor, mark=square*, mark size=1.25pt, line width=0.9pt, forget plot] coordinates {(100, 2.04000000e-16) (10000, 5.04000000e-16) (1000000, 5.59000000e-16) (100000000, 7.50000000e-16) (10000000000, 3.29000000e-16) (1000000000000, 4.12000000e-15) (100000000000000, 5.59000000e-16) (10000000000000000, 1.04000000e-15)};
\addplot+[color=LogColor, mark=triangle*, mark size=1.25pt, line width=0.9pt, forget plot] coordinates {(100, 2.73000000e-16) (10000, 1.03000000e-15) (1000000, 4.63000000e-16) (100000000, 1.52000000e-15) (10000000000, 1.13000000e-15) (1000000000000, 2.88000000e-13) (100000000000000, 1.61000000e-15) (10000000000000000, 9.68000000e-15)};
\addplot+[color=UniformColor, mark=diamond*, mark size=1.25pt, line width=0.9pt, forget plot] coordinates {(100, 1.54000000e-15) (10000, 5.28000000e-16) (1000000, 5.41000000e-16) (100000000, 6.81000000e-16) (10000000000, 1.89000000e-16) (1000000000000, 1.89000000e-16) (100000000000000, 3.53000000e-16) (10000000000000000, 5.90000000e-16)};
\addplot+[color=OneClusterColor, mark=pentagon*, mark size=1.25pt, line width=0.9pt, forget plot] coordinates {(100, 5.91000000e-16) (10000, 3.86000000e-16) (1000000, 1.36000000e-16) (100000000, 2.71000000e-16) (10000000000, 1.43000000e-16) (1000000000000, 2.84000000e-14) (100000000000000, 1.73000000e-13) (10000000000000000, 7.56000000e-14)};
\addplot+[color=TwoClusterColor, mark=o, mark size=1.25pt, line width=0.9pt, forget plot] coordinates {(100, 1.59000000e-15) (10000, 8.30000000e-16) (1000000, 3.39000000e-16) (100000000, 3.48000000e-16) (10000000000, 3.50000000e-16) (1000000000000, 3.55000000e-16) (100000000000000, 1.22000000e-15) (10000000000000000, 8.94000000e-15)};
\end{groupplot}
\node[anchor=south] at ([yshift=2mm]current bounding box.north) {\pgfplotslegendfromname{legendPSD1}};
\end{tikzpicture}
\caption{FT-Cholesky reconstruction residuals for positive semidefinite matrices with nullity one. Left: unblocked FT-Cholesky; right: blocked FT-Cholesky with block size $b=64$. The relative Frobenius residual is plotted against the target condition number. Each point represents one realization.}
\label{fig:ftchol-trend-psd1}
\end{figure}

Observe that the symmetric positive definite set of experiments produces reconstruction errors that are close to the numerical floor of Float64 representations (See \cref{fig:ftchol-trend-spd}). On the other hand, matrices with nullities have relatively worse reconstruction errors. The issue is particularly dire for clustered spectral shapes. This can be attributed to the large evolution times (here $T = 10^{14}$) used in these experiments. A modified algorithm that finds the Pareto front over $T$ for which the reconstruction errors are optimal can ameliorate some of these reconstruction issues. 

\begin{figure}[H]
\centering
\begin{tikzpicture}
\begin{groupplot}[
    group style={group size=2 by 1, horizontal sep=1.15cm},
    width=0.435\textwidth,
    height=0.31\textwidth,
    xmode=log,
    ymode=log,
    xmin=1e2, xmax=1e16,
    ymin=5e-17, ymax=1e-11,
    xtick={1e2,1e4,1e6,1e8,1e10,1e12,1e14,1e16},
    xticklabels={$10^2$,$10^4$,$10^6$,$10^8$,$10^{10}$,$10^{12}$,$10^{14}$,$10^{16}$},
    ytick={1e-16,1e-15,1e-14,1e-13,1e-12,1e-11},
    yticklabels={$10^{-16}$,$10^{-15}$,$10^{-14}$,$10^{-13}$,$10^{-12}$,$10^{-11}$},
    grid=major,
    major grid style={gray!25, line width=0.25pt},
    axis line style={gray!65},
    tick style={gray!65},
    tick label style={font=\small},
    label style={font=\small},
    title style={font=\small},
    xlabel={Condition number $\kappa$},
    scaled ticks=false,
    clip=true,
]
\nextgroupplot[
    title={Unblocked},
    ylabel={Relative Frobenius residual},
    legend to name=legendPSD2,
    legend columns=6,
    legend style={draw=none, fill=none, font=\small, /tikz/every even column/.append style={column sep=0.18cm}},
]
\addplot+[color=ExpColor, mark=*, mark size=1.25pt, line width=0.9pt] coordinates {(100, 7.75000000e-16) (10000, 1.57000000e-15) (1000000, 1.22000000e-15) (100000000, 3.74000000e-14) (10000000000, 1.24000000e-14) (1000000000000, 1.74000000e-15) (100000000000000, 2.69000000e-15) (10000000000000000, 1.04000000e-14)};
\addlegendentry{Exponential}
\addplot+[color=PolyColor, mark=square*, mark size=1.25pt, line width=0.9pt] coordinates {(100, 2.04000000e-16) (10000, 2.54000000e-15) (1000000, 8.90000000e-16) (100000000, 9.11000000e-16) (10000000000, 8.70000000e-16) (1000000000000, 7.78000000e-15) (100000000000000, 1.23000000e-15) (10000000000000000, 1.09000000e-15)};
\addlegendentry{Polynomial}
\addplot+[color=LogColor, mark=triangle*, mark size=1.25pt, line width=0.9pt] coordinates {(100, 2.28000000e-15) (10000, 1.11000000e-14) (1000000, 1.74000000e-15) (100000000, 1.29000000e-15) (10000000000, 2.07000000e-15) (1000000000000, 1.79000000e-15) (100000000000000, 1.70000000e-15) (10000000000000000, 9.97000000e-15)};
\addlegendentry{Log-uniform}
\addplot+[color=UniformColor, mark=diamond*, mark size=1.25pt, line width=0.9pt] coordinates {(100, 7.38000000e-16) (10000, 5.41000000e-16) (1000000, 7.10000000e-16) (100000000, 1.98000000e-14) (10000000000, 1.16000000e-15) (1000000000000, 9.56000000e-16) (100000000000000, 1.02000000e-15) (10000000000000000, 1.14000000e-14)};
\addlegendentry{Uniform}
\addplot+[color=OneClusterColor, mark=pentagon*, mark size=1.25pt, line width=0.9pt] coordinates {(100, 3.20000000e-16) (10000, 6.89000000e-16) (1000000, 2.13000000e-15) (100000000, 1.80000000e-15) (10000000000, 1.47000000e-15) (1000000000000, 4.76000000e-16) (100000000000000, 7.36000000e-16) (10000000000000000, 1.18000000e-14)};
\addlegendentry{One-cluster}
\addplot+[color=TwoClusterColor, mark=o, mark size=1.25pt, line width=0.9pt] coordinates {(100, 7.34000000e-16) (10000, 5.30000000e-16) (1000000, 1.87000000e-15) (100000000, 2.30000000e-12) (10000000000, 4.26000000e-14) (1000000000000, 1.07000000e-15) (100000000000000, 1.28000000e-15) (10000000000000000, 9.16000000e-15)};
\addlegendentry{Two-cluster}
\nextgroupplot[
    title={Blocked ($b=64$)},
    ylabel={},
]
\addplot+[color=ExpColor, mark=*, mark size=1.25pt, line width=0.9pt, forget plot] coordinates {(100, 7.76000000e-16) (10000, 1.24000000e-15) (1000000, 1.31000000e-15) (100000000, 1.26000000e-15) (10000000000, 7.82000000e-15) (1000000000000, 1.73000000e-15) (100000000000000, 2.42000000e-15) (10000000000000000, 1.04000000e-14)};
\addplot+[color=PolyColor, mark=square*, mark size=1.25pt, line width=0.9pt, forget plot] coordinates {(100, 1.80000000e-16) (10000, 3.94000000e-14) (1000000, 9.58000000e-16) (100000000, 8.77000000e-16) (10000000000, 6.79000000e-16) (1000000000000, 1.31000000e-15) (100000000000000, 1.22000000e-15) (10000000000000000, 1.17000000e-15)};
\addplot+[color=LogColor, mark=triangle*, mark size=1.25pt, line width=0.9pt, forget plot] coordinates {(100, 6.63000000e-16) (10000, 1.79000000e-14) (1000000, 1.63000000e-15) (100000000, 1.30000000e-15) (10000000000, 1.93000000e-15) (1000000000000, 2.02000000e-15) (100000000000000, 1.78000000e-15) (10000000000000000, 9.96000000e-15)};
\addplot+[color=UniformColor, mark=diamond*, mark size=1.25pt, line width=0.9pt, forget plot] coordinates {(100, 6.92000000e-16) (10000, 4.24000000e-16) (1000000, 2.89000000e-15) (100000000, 1.30000000e-14) (10000000000, 4.99000000e-15) (1000000000000, 1.29000000e-15) (100000000000000, 8.56000000e-16) (10000000000000000, 9.10000000e-15)};
\addplot+[color=OneClusterColor, mark=pentagon*, mark size=1.25pt, line width=0.9pt, forget plot] coordinates {(100, 3.68000000e-16) (10000, 2.81000000e-15) (1000000, 1.99000000e-16) (100000000, 1.23000000e-15) (10000000000, 5.94000000e-15) (1000000000000, 3.10000000e-16) (100000000000000, 5.93000000e-16) (10000000000000000, 1.02000000e-14)};
\addplot+[color=TwoClusterColor, mark=o, mark size=1.25pt, line width=0.9pt, forget plot] coordinates {(100, 7.00000000e-16) (10000, 5.05000000e-16) (1000000, 5.97000000e-16) (100000000, 8.31000000e-16) (10000000000, 7.37000000e-14) (1000000000000, 8.32000000e-16) (100000000000000, 1.26000000e-15) (10000000000000000, 9.20000000e-15)};
\end{groupplot}
\node[anchor=south] at ([yshift=2mm]current bounding box.north) {\pgfplotslegendfromname{legendPSD2}};
\end{tikzpicture}
\caption{FT-Cholesky reconstruction residuals for positive semidefinite matrices with nullity two. Left: unblocked FT-Cholesky; right: blocked FT-Cholesky with block size $b=64$. The relative Frobenius residual is plotted against the target condition number. Each point represents one realization.}
\label{fig:ftchol-trend-psd2}
\end{figure}

Furthermore, some of these matrices are numerically indefinite in Float64 representations. This will consequently spawn negative pivots in the FT-Cholesky algorithm -- a regime where the FT-Cholesky map is not well-posed. Nonetheless, when the negative pivots are in a sufficiently small $\epsilon$ ball around zero, the FT-Cholesky algorithm works without introducing spurious compounding errors. 

\begin{figure}[H]
\centering
\begin{tikzpicture}
\begin{groupplot}[
    group style={group size=2 by 1, horizontal sep=1.15cm},
    width=0.435\textwidth,
    height=0.31\textwidth,
    xmode=log,
    ymode=log,
    xmin=1e2, xmax=1e16,
    ymin=5e-17, ymax=1e-11,
    xtick={1e2,1e4,1e6,1e8,1e10,1e12,1e14,1e16},
    xticklabels={$10^2$,$10^4$,$10^6$,$10^8$,$10^{10}$,$10^{12}$,$10^{14}$,$10^{16}$},
    ytick={1e-16,1e-15,1e-14,1e-13,1e-12,1e-11},
    yticklabels={$10^{-16}$,$10^{-15}$,$10^{-14}$,$10^{-13}$,$10^{-12}$,$10^{-11}$},
    grid=major,
    major grid style={gray!25, line width=0.25pt},
    axis line style={gray!65},
    tick style={gray!65},
    tick label style={font=\small},
    label style={font=\small},
    title style={font=\small},
    xlabel={Condition number $\kappa$},
    scaled ticks=false,
    clip=true,
]
\nextgroupplot[
    title={Unblocked},
    ylabel={Relative Frobenius residual},
    legend to name=legendPSD5,
    legend columns=6,
    legend style={draw=none, fill=none, font=\small, /tikz/every even column/.append style={column sep=0.18cm}},
]
\addplot+[color=ExpColor, mark=*, mark size=1.25pt, line width=0.9pt] coordinates {(100, 1.08000000e-14) (10000, 3.56000000e-15) (1000000, 4.99000000e-15) (100000000, 3.18000000e-15) (10000000000, 1.06000000e-14) (1000000000000, 3.90000000e-15) (100000000000000, 5.75000000e-15) (10000000000000000, 1.14000000e-14)};
\addlegendentry{Exponential}
\addplot+[color=PolyColor, mark=square*, mark size=1.25pt, line width=0.9pt] coordinates {(100, 1.17000000e-15) (10000, 2.36000000e-12) (1000000, 1.08000000e-15) (100000000, 1.86000000e-15) (10000000000, 1.24000000e-13) (1000000000000, 8.13000000e-14) (100000000000000, 1.87000000e-15) (10000000000000000, 2.32000000e-15)};
\addlegendentry{Polynomial}
\addplot+[color=LogColor, mark=triangle*, mark size=1.25pt, line width=0.9pt] coordinates {(100, 1.86000000e-15) (10000, 2.07000000e-15) (1000000, 2.53000000e-15) (100000000, 5.97000000e-15) (10000000000, 4.85000000e-15) (1000000000000, 3.61000000e-15) (100000000000000, 4.32000000e-15) (10000000000000000, 1.18000000e-14)};
\addlegendentry{Log-uniform}
\addplot+[color=UniformColor, mark=diamond*, mark size=1.25pt, line width=0.9pt] coordinates {(100, 1.99000000e-14) (10000, 6.30000000e-15) (1000000, 7.01000000e-15) (100000000, 1.38000000e-15) (10000000000, 1.34000000e-15) (1000000000000, 1.68000000e-13) (100000000000000, 3.54000000e-15) (10000000000000000, 1.48000000e-15)};
\addlegendentry{Uniform}
\addplot+[color=OneClusterColor, mark=pentagon*, mark size=1.25pt, line width=0.9pt] coordinates {(100, 9.99000000e-16) (10000, 1.66000000e-12) (1000000, 4.94000000e-15) (100000000, 4.81000000e-15) (10000000000, 2.83000000e-15) (1000000000000, 6.47000000e-15) (100000000000000, 2.58000000e-15) (10000000000000000, 4.45000000e-12)};
\addlegendentry{One-cluster}
\addplot+[color=TwoClusterColor, mark=o, mark size=1.25pt, line width=0.9pt] coordinates {(100, 8.42000000e-16) (10000, 8.25000000e-14) (1000000, 9.05000000e-16) (100000000, 2.16000000e-14) (10000000000, 1.26000000e-15) (1000000000000, 5.94000000e-14) (100000000000000, 1.44000000e-15) (10000000000000000, 9.48000000e-15)};
\addlegendentry{Two-cluster}
\nextgroupplot[
    title={Blocked ($b=64$)},
    ylabel={},
]
\addplot+[color=ExpColor, mark=*, mark size=1.25pt, line width=0.9pt, forget plot] coordinates {(100, 6.09000000e-15) (10000, 5.97000000e-15) (1000000, 4.28000000e-15) (100000000, 4.11000000e-15) (10000000000, 1.10000000e-14) (1000000000000, 5.62000000e-15) (100000000000000, 6.32000000e-15) (10000000000000000, 1.14000000e-14)};
\addplot+[color=PolyColor, mark=square*, mark size=1.25pt, line width=0.9pt, forget plot] coordinates {(100, 1.54000000e-15) (10000, 7.99000000e-16) (1000000, 1.42000000e-15) (100000000, 1.70000000e-15) (10000000000, 1.60000000e-13) (1000000000000, 2.52000000e-15) (100000000000000, 1.74000000e-15) (10000000000000000, 2.34000000e-15)};
\addplot+[color=LogColor, mark=triangle*, mark size=1.25pt, line width=0.9pt, forget plot] coordinates {(100, 6.63000000e-15) (10000, 2.08000000e-15) (1000000, 2.41000000e-15) (100000000, 4.27000000e-15) (10000000000, 5.55000000e-15) (1000000000000, 3.22000000e-15) (100000000000000, 4.30000000e-15) (10000000000000000, 1.18000000e-14)};
\addplot+[color=UniformColor, mark=diamond*, mark size=1.25pt, line width=0.9pt, forget plot] coordinates {(100, 1.98000000e-14) (10000, 6.10000000e-15) (1000000, 2.24000000e-15) (100000000, 9.78000000e-16) (10000000000, 9.24000000e-16) (1000000000000, 1.71000000e-13) (100000000000000, 3.18000000e-15) (10000000000000000, 4.46000000e-15)};
\addplot+[color=OneClusterColor, mark=pentagon*, mark size=1.25pt, line width=0.9pt, forget plot] coordinates {(100, 8.54000000e-16) (10000, 1.96000000e-12) (1000000, 5.53000000e-15) (100000000, 2.42000000e-15) (10000000000, 2.30000000e-15) (1000000000000, 3.58000000e-15) (100000000000000, 2.03000000e-15) (10000000000000000, 5.80000000e-12)};
\addplot+[color=TwoClusterColor, mark=o, mark size=1.25pt, line width=0.9pt, forget plot] coordinates {(100, 7.23000000e-16) (10000, 1.88000000e-14) (1000000, 1.05000000e-15) (100000000, 1.57000000e-14) (10000000000, 1.28000000e-15) (1000000000000, 8.02000000e-15) (100000000000000, 1.43000000e-15) (10000000000000000, 9.61000000e-15)};
\end{groupplot}
\node[anchor=south] at ([yshift=2mm]current bounding box.north) {\pgfplotslegendfromname{legendPSD5}};
\end{tikzpicture}
\caption{FT-Cholesky reconstruction residuals for positive semidefinite matrices with nullity five. Left: unblocked FT-Cholesky; right: blocked FT-Cholesky with block size $b=64$. The relative Frobenius residual is plotted against the target condition number. Each point represents one realization.}
\label{fig:ftchol-trend-psd5}
\end{figure}

\begin{figure}[H]
\centering
\begin{tikzpicture}
\begin{groupplot}[
    group style={group size=2 by 1, horizontal sep=1.15cm},
    width=0.435\textwidth,
    height=0.31\textwidth,
    xmode=log,
    ymode=log,
    xmin=1e2, xmax=1e16,
    ymin=5e-17, ymax=1e-11,
    xtick={1e2,1e4,1e6,1e8,1e10,1e12,1e14,1e16},
    xticklabels={$10^2$,$10^4$,$10^6$,$10^8$,$10^{10}$,$10^{12}$,$10^{14}$,$10^{16}$},
    ytick={1e-16,1e-15,1e-14,1e-13,1e-12,1e-11},
    yticklabels={$10^{-16}$,$10^{-15}$,$10^{-14}$,$10^{-13}$,$10^{-12}$,$10^{-11}$},
    grid=major,
    major grid style={gray!25, line width=0.25pt},
    axis line style={gray!65},
    tick style={gray!65},
    tick label style={font=\small},
    label style={font=\small},
    title style={font=\small},
    xlabel={Condition number $\kappa$},
    scaled ticks=false,
    clip=true,
]
\nextgroupplot[
    title={Unblocked},
    ylabel={Relative Frobenius residual},
    legend to name=legendPSD10,
    legend columns=6,
    legend style={draw=none, fill=none, font=\small, /tikz/every even column/.append style={column sep=0.18cm}},
]
\addplot+[color=ExpColor, mark=*, mark size=1.25pt, line width=0.9pt] coordinates {(100, 3.51000000e-15) (10000, 9.46000000e-15) (1000000, 3.60000000e-14) (100000000, 5.00000000e-15) (10000000000, 7.03000000e-15) (1000000000000, 5.16000000e-15) (100000000000000, 1.13000000e-14) (10000000000000000, 1.33000000e-14)};
\addlegendentry{Exponential}
\addplot+[color=PolyColor, mark=square*, mark size=1.25pt, line width=0.9pt] coordinates {(100, 1.18000000e-15) (10000, 6.10000000e-15) (1000000, 1.59000000e-13) (100000000, 2.92000000e-15) (10000000000, 7.02000000e-15) (1000000000000, 3.13000000e-15) (100000000000000, 1.90000000e-13) (10000000000000000, 2.78000000e-15)};
\addlegendentry{Polynomial}
\addplot+[color=LogColor, mark=triangle*, mark size=1.25pt, line width=0.9pt] coordinates {(100, 2.57000000e-15) (10000, 2.33000000e-14) (1000000, 5.93000000e-15) (100000000, 4.40000000e-15) (10000000000, 1.55000000e-12) (1000000000000, 9.14000000e-15) (100000000000000, 7.51000000e-15) (10000000000000000, 1.33000000e-14)};
\addlegendentry{Log-uniform}
\addplot+[color=UniformColor, mark=diamond*, mark size=1.25pt, line width=0.9pt] coordinates {(100, 5.26000000e-15) (10000, 6.27000000e-15) (1000000, 1.45000000e-14) (100000000, 8.27000000e-15) (10000000000, 2.70000000e-14) (1000000000000, 1.43000000e-15) (100000000000000, 1.50000000e-15) (10000000000000000, 2.76000000e-15)};
\addlegendentry{Uniform}
\addplot+[color=OneClusterColor, mark=pentagon*, mark size=1.25pt, line width=0.9pt] coordinates {(100, 5.14000000e-14) (10000, 6.26000000e-14) (1000000, 7.71000000e-15) (100000000, 1.73000000e-15) (10000000000, 1.91000000e-15) (1000000000000, 2.21000000e-14) (100000000000000, 5.87000000e-15) (10000000000000000, 5.23000000e-14)};
\addlegendentry{One-cluster}
\addplot+[color=TwoClusterColor, mark=o, mark size=1.25pt, line width=0.9pt] coordinates {(100, 3.59000000e-15) (10000, 2.57000000e-14) (1000000, 3.99000000e-15) (100000000, 5.90000000e-15) (10000000000, 2.00000000e-13) (1000000000000, 2.46000000e-15) (100000000000000, 1.72000000e-15) (10000000000000000, 1.16000000e-14)};
\addlegendentry{Two-cluster}
\nextgroupplot[
    title={Blocked ($b=64$)},
    ylabel={},
]
\addplot+[color=ExpColor, mark=*, mark size=1.25pt, line width=0.9pt, forget plot] coordinates {(100, 2.75000000e-15) (10000, 2.68000000e-14) (1000000, 6.09000000e-15) (100000000, 4.83000000e-15) (10000000000, 6.79000000e-15) (1000000000000, 5.46000000e-15) (100000000000000, 1.08000000e-14) (10000000000000000, 1.33000000e-14)};
\addplot+[color=PolyColor, mark=square*, mark size=1.25pt, line width=0.9pt, forget plot] coordinates {(100, 1.35000000e-15) (10000, 5.68000000e-15) (1000000, 2.54000000e-13) (100000000, 2.78000000e-15) (10000000000, 1.63000000e-14) (1000000000000, 4.15000000e-15) (100000000000000, 4.44000000e-15) (10000000000000000, 3.08000000e-15)};
\addplot+[color=LogColor, mark=triangle*, mark size=1.25pt, line width=0.9pt, forget plot] coordinates {(100, 2.05000000e-15) (10000, 7.30000000e-14) (1000000, 4.82000000e-15) (100000000, 5.89000000e-15) (10000000000, 2.41000000e-13) (1000000000000, 9.15000000e-14) (100000000000000, 8.81000000e-15) (10000000000000000, 1.33000000e-14)};
\addplot+[color=UniformColor, mark=diamond*, mark size=1.25pt, line width=0.9pt, forget plot] coordinates {(100, 2.11000000e-14) (10000, 2.48000000e-15) (1000000, 1.16000000e-14) (100000000, 5.29000000e-15) (10000000000, 1.16000000e-14) (1000000000000, 1.40000000e-15) (100000000000000, 1.62000000e-15) (10000000000000000, 2.80000000e-15)};
\addplot+[color=OneClusterColor, mark=pentagon*, mark size=1.25pt, line width=0.9pt, forget plot] coordinates {(100, 4.42000000e-14) (10000, 5.66000000e-14) (1000000, 5.75000000e-15) (100000000, 1.91000000e-15) (10000000000, 1.66000000e-15) (1000000000000, 2.11000000e-14) (100000000000000, 5.79000000e-15) (10000000000000000, 4.72000000e-14)};
\addplot+[color=TwoClusterColor, mark=o, mark size=1.25pt, line width=0.9pt, forget plot] coordinates {(100, 3.84000000e-15) (10000, 1.94000000e-14) (1000000, 6.36000000e-15) (100000000, 3.11000000e-15) (10000000000, 1.26000000e-13) (1000000000000, 8.31000000e-15) (100000000000000, 1.71000000e-15) (10000000000000000, 1.22000000e-14)};
\end{groupplot}
\node[anchor=south] at ([yshift=2mm]current bounding box.north) {\pgfplotslegendfromname{legendPSD10}};
\end{tikzpicture}
\caption{FT-Cholesky reconstruction residuals for positive semidefinite matrices with nullity ten. Left: unblocked FT-Cholesky; right: blocked FT-Cholesky with block size $b=64$. The relative Frobenius residual is plotted against the target condition number. Each point represents one realization.}
\label{fig:ftchol-trend-psd10}
\end{figure}

\subsection{Applications}
\label{sec:applications}
In the rest of this section, we will test the FT-Cholesky algorithms on matrices that arise in different applications. Particular attention is given to a litany of kernel matrices and correspondingly to finite-difference matrices of linear elliptic differential operators. We are again interested in the relative Frobenius residual of the factorization. Comparisons are made with \texttt{xPOTRF}, \texttt{xPSTRF} and its blocked implementations. 

\subsubsection{Kernel matrices}
\label{subsec:kernel-applications}
Kernel estimation \cite{scholkopf2002learning}, RBF interpolation \cite{schaback2006kernel}, and Gaussian process regression \cite{rasmussen2006gaussian} are related and, under some assumptions, equivalent methods for estimating a function from noisy observations. Let $\mathcal{X} \subset \mathbb{R}^d$ and $\mathcal{Y} \subset \mathbb{R}$ denote the input and output spaces, respectively. A target function $f : \mathcal{X} \mapsto \mathcal{Y}$ is observed on certain collocation points $\hat{\mathcal{X}} \subset \mathcal{X}$, denoted by $\hat{\mathcal{Y}} := f(\hat{\mathcal{X}}) + \mathcal{E}$ where $\mathcal{E}$ denotes the appropriate pointwise i.i.d. Gaussian noise with variance $\lambda$. The kernel estimation problem tries to approximate $f$ using an estimator:
\begin{align}
    \begin{gathered}
        g(\cdot) := \sum_{i=1}^M a_i k(\cdot, \hat{\mathcal{X}_i}). 
    \end{gathered}
\end{align}
Here, $k : \mathcal{X} \times \mathcal{X} \mapsto \mathbb{R}$ is a positive-semi-definite function known as the kernel function. Under certain assumptions, these kernel functions are often Green's functions of certain linear elliptic differential operators. Common examples of such kernel functions include the Matern kernel, the Wendland kernel, the Gaussian kernel \cite{wendland2005scattered} and the like. For a comprehensive review on kernel functions, see \cite{wendland2005scattered}. To estimate $g$, one solves the following linear least squares problem:
\begin{align}
    \begin{gathered}
        \arg \min_{a \in \mathbb{R}^M}  \|k(\hat{\mathcal{X}}, \hat{\mathcal{X}}) a - f(\hat{\mathcal{X}})\|^2_2.
    \end{gathered}
\end{align}

Furthermore, the average incurred in estimating  $g$ with finite information is given by its power function (it's also the posterior covariance in Gaussian process regression). For $x,y \in \mathcal{X}$, it is defined as:
\begin{align}
    \begin{gathered}
        \mathcal{P}(x,y) := k(x,y) - k(x, \hat{\mathcal{X}}) (k(\hat{\mathcal{X}},\hat{\mathcal{X}})+ \lambda I)^{-1} k(\hat{\mathcal{X}}, y).  
    \end{gathered}
\end{align}

Both these quantities -- $a, \mathcal{P}$ can be easily estimated once the Cholesky factor of $k(\hat{\mathcal{X}},\hat{\mathcal{X}})$ is estimated. However, in applications where $\lambda$ is relatively small, the canonical Cholesky implementation terminates prematurely. Therefore, we will carry out this factorization using our finite-time algorithms. 

\begin{table}[htbp]
\centering
\caption{Reconstruction errors for singular kernel matrices generated by exact repeated points. All matrices have size \(n=500\). Revised unblocked FT-Cholesky uses \(T=10^{14}\). The error is \(\|K-LL^\top\|_F/\|K\|_F\), accounting for permutation for pivoted factors. Standard and pivoted Cholesky use \texttt{DPOTF2} and \texttt{DPSTF2}; pivoted rank uses LAPACK's default tolerance.}
\label{tab:kernel-reconstruction-unblocked}
\begin{tabular}{lrrlrr}
\toprule
Kernel & Nullity & Finite-time error & \texttt{DPOTF2} & Pivoted rank & Pivoted error \\
\midrule
Mat\'ern-$5/2$ & 1 & $3.629\times10^{-16}$ & Pass & 499 & $3.241\times10^{-16}$ \\
 & 2 & $3.596\times10^{-16}$ & Fail & 498 & $3.253\times10^{-16}$ \\
 & 5 & $3.643\times10^{-16}$ & Fail & 495 & $3.216\times10^{-16}$ \\
 & 10 & $3.664\times10^{-16}$ & Fail & 490 & $3.253\times10^{-16}$ \\
\midrule
IWP-1 & 1 & $5.434\times10^{-16}$ & Pass & 499 & $3.073\times10^{-16}$ \\
 & 2 & $5.640\times10^{-16}$ & Fail & 498 & $3.052\times10^{-16}$ \\
 & 5 & $5.955\times10^{-16}$ & Fail & 495 & $3.054\times10^{-16}$ \\
 & 10 & $6.592\times10^{-16}$ & Fail & 490 & $3.082\times10^{-16}$ \\
\midrule
Wendland & 1 & $1.764\times10^{-16}$ & Fail & 499 & $3.515\times10^{-16}$ \\
 & 2 & $1.824\times10^{-16}$ & Fail & 498 & $3.500\times10^{-16}$ \\
 & 5 & $1.949\times10^{-16}$ & Fail & 495 & $3.514\times10^{-16}$ \\
 & 10 & $2.229\times10^{-16}$ & Fail & 490 & $3.552\times10^{-16}$ \\
\bottomrule
\end{tabular}
\end{table}

The reconstruction errors of kernel matrices due to Matern, integrated-Wiener process \cite{oksendal2003stochastic} and Wendland priors using the finite-time Cholesky algorithm and its blocked version are furnished in \cref{tab:kernel-reconstruction-unblocked} and \cref{tab:kernel-reconstruction-blocked}, respectively. These matrices are constructed using randomly sampled points on an interval of $\mathbb{R}$. The nullities of these matrices are controlled by overlapping points in the collocation set. While this appears contrived, uniform random sampling of collocation points has similar tendencies. 

\begin{table}[htbp]
\centering
\caption{Reconstruction errors for singular kernel matrices generated by exact repeated points. All matrices have size \(n=500\). Revised blocked FT-Cholesky uses \(T=10^{14}\) and block size \(b=64\). The error is \(\|K-LL^\top\|_F/\|K\|_F\), accounting for permutation for pivoted factors. Standard and pivoted Cholesky use \texttt{DPOTRF} and \texttt{DPSTRF}; pivoted rank uses LAPACK's default tolerance.}
\label{tab:kernel-reconstruction-blocked}
\begin{tabular}{lrrlrr}
\toprule
Kernel & Nullity & Finite-time error & \texttt{DPOTRF} & Pivoted rank & Pivoted error \\
\midrule
Mat\'ern-$5/2$ & 1 & $2.284\times10^{-16}$ & Pass & 499 & $2.113\times10^{-16}$ \\
 & 2 & $2.301\times10^{-16}$ & Fail & 498 & $2.140\times10^{-16}$ \\
 & 5 & $2.343\times10^{-16}$ & Fail & 495 & $2.129\times10^{-16}$ \\
 & 10 & $2.400\times10^{-16}$ & Fail & 490 & $2.142\times10^{-16}$ \\
\midrule
IWP-1 & 1 & $5.177\times10^{-16}$ & Pass & 499 & $1.834\times10^{-16}$ \\
 & 2 & $5.308\times10^{-16}$ & Fail & 498 & $1.821\times10^{-16}$ \\
 & 5 & $5.715\times10^{-16}$ & Fail & 495 & $1.817\times10^{-16}$ \\
 & 10 & $6.375\times10^{-16}$ & Fail & 490 & $1.833\times10^{-16}$ \\
\midrule
Wendland & 1 & $1.567\times10^{-16}$ & Fail & 499 & $2.791\times10^{-16}$ \\
 & 2 & $1.657\times10^{-16}$ & Fail & 498 & $2.825\times10^{-16}$ \\
 & 5 & $1.819\times10^{-16}$ & Fail & 495 & $2.801\times10^{-16}$ \\
 & 10 & $2.109\times10^{-16}$ & Fail & 490 & $2.799\times10^{-16}$ \\
\bottomrule
\end{tabular}
\end{table}

Note that both of our methods are competitive with pivoted routines while preserving the memory access patterns of the canonical Cholesky method. While pivoted routines have the advantage of being rank-revealing, it is seldom interesting in practice to estimate the rank of these Gram matrices. 

\begin{remark}
Approximation bounds (along with empirical investigations) on the estimating  $a, \mathcal{P}$ will be pursued in a subsequent work. 
\end{remark}

\subsubsection{Hermite kernel matrices}
\label{subsec:hermite-kernels}

A natural generalization of kernel estimation involves approximation, where $f : \mathcal{X} \mapsto \mathcal{Y}$ is smooth, i.e. $C^k(\mathcal{X}, \mathcal{Y})$. Furthermore, we assume that, in addition to function values on $\hat{\mathcal{X}}$, we also measure certain derivatives of $f$, namely $D f, D^2 f, \dots$ at the same points. We denote the measurement sets by $\hat{\mathcal{Y}}_{k=0}^s := D^k f(\hat{\mathcal{X}}) + \mathcal{E}_k$ where $\mathcal{E}_k$ is again pointwise i.i.d Gaussian noise with variances $\lambda_k > 0$. The estimator is specified using the following ansatz:
\begin{align}
    \begin{gathered}
        h(x) := \sum_{k=0}^s \sum_{j=1}^M b_j^k \; [D^k]_{1} k(x, \mathcal{X}_j). 
    \end{gathered}
\end{align}
Here, $[D^k]_{1}$ denotes the $k^{th}$ derivative operator with respect to the first argument of the kernel $k$, and $b_{j}^k \in \mathbb{R}$ denotes the coefficients to be estimated. To estimate these parameters, we consider the derivatives of this ansatz:
\begin{align}
    \begin{gathered}
        [D^k]_{2} h(x) = \sum_{k=0}^s \sum_{j=1}^M b_j^k \; [D^k]_{2}\;  [D^k]_{1} k(x, \mathcal{X}_j).
    \end{gathered}
\end{align}
Taken together, this estimation procedure with derivative information induces a Gram matrix, which we denote by:
\begin{align}
    \Theta_{\hat{\mathcal{X}},s,M}
    := \begin{bmatrix}
        k(\hat{\mathcal{X}},\hat{\mathcal{X}})
        & [D^1]_2 k(\hat{\mathcal{X}},\hat{\mathcal{X}})
        & \cdots
        & [D^s]_2 k(\hat{\mathcal{X}},\hat{\mathcal{X}}) \\
        [D^1]_1 k(\hat{\mathcal{X}},\hat{\mathcal{X}})
        & [D^1]_1[D^1]_2 k(\hat{\mathcal{X}},\hat{\mathcal{X}})
        & \cdots
        & [D^1]_1[D^s]_2 k(\hat{\mathcal{X}},\hat{\mathcal{X}}) \\
        \vdots & \vdots & \ddots & \vdots \\
        [D^s]_1 k(\hat{\mathcal{X}},\hat{\mathcal{X}})
        & [D^s]_1[D^1]_2 k(\hat{\mathcal{X}},\hat{\mathcal{X}})
        & \cdots
        & [D^s]_1[D^s]_2 k(\hat{\mathcal{X}},\hat{\mathcal{X}})
    \end{bmatrix}.
\end{align}
Here, each block say $k(\hat{\mathcal{X}},\hat{\mathcal{X}}) \in \mathbb{R}^{M \times M}$. Furthermore, each block is not positive semi-definite. On the other hand, $\Theta_{\hat{\mathcal{X}}, s, M} \succeq 0$. We will refer to such structures as Hermite-kernel matrices, eponymous with the Hermite-Birkhoff interpolation problem. An immediate application of this is in solving partial differential equations using collocation \cite{chen2021solving}. For instance, when solving the Poisson equation:
\begin{align}
    \begin{gathered}
        -\Delta u(x) = p(x), \; x \in \mathcal{X}, \\ 
        u(x) = 0, \; \lim \|x\| \mapsto \infty. 
    \end{gathered}
\end{align}
where $p \in L^2(\Omega)$, one assembles the following Hermite-kernel matrix:
\begin{align}
    \begin{gathered}
        \Theta_{\Delta} :=
        \begin{bmatrix}
            k(\hat{\mathcal{X}},\hat{\mathcal{X}})
            & \Delta_2 k(\hat{\mathcal{X}},\hat{\mathcal{X}}) 
            \\ 
            \Delta_1 k(\hat{\mathcal{X}},\hat{\mathcal{X}}) & \Delta_1 \Delta_2 k(\hat{\mathcal{X}},\hat{\mathcal{X}}) 
        \end{bmatrix}
    \end{gathered}
\end{align}
Depending on the smoothness of the kernel, these matrices can be highly ill-conditioned. In the following, we will compare the reconstruction accuracies corresponding to different integer Matern kernels. To control the null space of these matrices, we deliberately introduce duplicate points when the collocation points are sampled. The relative reconstruction accuracies of the unblocked and blocked finite-time Cholesky implementations are furnished in \cref{tab:poisson-hermite-unblocked} and \cref{tab:elliptic-reconstruction-blocked}, respectively. 

We observe that both the unblocked and blocked implementations successfully complete factorization for all the considered kernel matrices. Furthermore, it is competitive with LAPACK's pivoted Cholesky implementations. The un-pivoted implementations on the other hand, terminate early upon encountering a zero pivot. 

\begin{table}[htbp]
\centering
\caption{Reconstruction errors for mixed first-order Poisson Hermite Gram matrices with exact repeated operator observations. All matrices have size \(n=1000\). Revised unblocked FT-Cholesky uses \(T=10^{14}\). The error is \(\|K-LL^\top\|_F/\|K\|_F\), accounting for permutation for pivoted factors. Standard and pivoted Cholesky use \texttt{DPOTF2} and \texttt{DPSTF2}; pivoted rank uses LAPACK's default tolerance.}
\label{tab:poisson-hermite-unblocked}
\begin{tabular}{lrrlrr}
\toprule
Kernel & Repeats & Finite-time error & \texttt{DPOTF2} & Pivoted rank & Pivoted error \\
\midrule
Mat\'ern-$3/2$ & 0 & $5.157\times10^{-16}$ & Pass & 1000 & $5.955\times10^{-16}$ \\
 & 1 & $5.166\times10^{-16}$ & Fail & 999 & $5.875\times10^{-16}$ \\
 & 2 & $5.195\times10^{-16}$ & Fail & 998 & $5.954\times10^{-16}$ \\
 & 5 & $5.169\times10^{-16}$ & Fail & 995 & $5.886\times10^{-16}$ \\
 & 10 & $5.248\times10^{-16}$ & Fail & 990 & $5.847\times10^{-16}$ \\
\midrule
Mat\'ern-$5/2$ & 0 & $5.499\times10^{-16}$ & Pass & 1000 & $7.752\times10^{-16}$ \\
 & 1 & $5.520\times10^{-16}$ & Fail & 999 & $7.582\times10^{-16}$ \\
 & 2 & $5.518\times10^{-16}$ & Fail & 998 & $7.618\times10^{-16}$ \\
 & 5 & $5.483\times10^{-16}$ & Fail & 995 & $7.592\times10^{-16}$ \\
 & 10 & $5.480\times10^{-16}$ & Fail & 990 & $7.539\times10^{-16}$ \\
\midrule
Mat\'ern-$7/2$ & 0 & $5.539\times10^{-16}$ & Pass & 1000 & $6.322\times10^{-16}$ \\
 & 1 & $5.535\times10^{-16}$ & Fail & 999 & $6.304\times10^{-16}$ \\
 & 2 & $5.561\times10^{-16}$ & Pass & 998 & $6.320\times10^{-16}$ \\
 & 5 & $5.529\times10^{-16}$ & Fail & 995 & $6.281\times10^{-16}$ \\
 & 10 & $5.466\times10^{-16}$ & Fail & 990 & $6.324\times10^{-16}$ \\
\bottomrule
\end{tabular}
\end{table}

\begin{table}[htbp]
\centering
\caption{Reconstruction errors for mixed first-order Poisson Hermite Gram matrices with exact repeated operator observations. All matrices have size \(n=1000\). Revised blocked FT-Cholesky uses \(T=10^{14}\) and block size \(b=64\). The error is \(\|K-LL^\top\|_F/\|K\|_F\), accounting for permutation for pivoted factors. Standard and pivoted Cholesky use \texttt{DPOTRF} and \texttt{DPSTRF}; pivoted rank uses LAPACK's default tolerance.}
\label{tab:poisson-hermite-blocked}
\begin{tabular}{lrrlrr}
\toprule
Kernel & Repeats & Finite-time error & \texttt{DPOTRF} & Pivoted rank & Pivoted error \\
\midrule
Mat\'ern-$3/2$ & 0 & $2.068\times10^{-16}$ & Pass & 1000 & $5.113\times10^{-16}$ \\
 & 1 & $2.058\times10^{-16}$ & Fail & 999 & $5.091\times10^{-16}$ \\
 & 2 & $2.040\times10^{-16}$ & Fail & 998 & $5.118\times10^{-16}$ \\
 & 5 & $2.064\times10^{-16}$ & Fail & 995 & $5.099\times10^{-16}$ \\
 & 10 & $2.073\times10^{-16}$ & Fail & 990 & $5.027\times10^{-16}$ \\
\midrule
Mat\'ern-$5/2$ & 0 & $2.387\times10^{-16}$ & Pass & 1000 & $6.361\times10^{-16}$ \\
 & 1 & $2.384\times10^{-16}$ & Pass & 999 & $6.239\times10^{-16}$ \\
 & 2 & $2.386\times10^{-16}$ & Fail & 998 & $6.241\times10^{-16}$ \\
 & 5 & $2.390\times10^{-16}$ & Fail & 995 & $6.258\times10^{-16}$ \\
 & 10 & $2.403\times10^{-16}$ & Fail & 990 & $6.181\times10^{-16}$ \\
\midrule
Mat\'ern-$7/2$ & 0 & $2.599\times10^{-16}$ & Pass & 1000 & $5.502\times10^{-16}$ \\
 & 1 & $2.593\times10^{-16}$ & Fail & 999 & $5.489\times10^{-16}$ \\
 & 2 & $2.643\times10^{-16}$ & Pass & 998 & $5.513\times10^{-16}$ \\
 & 5 & $2.615\times10^{-16}$ & Fail & 995 & $5.466\times10^{-16}$ \\
 & 10 & $2.607\times10^{-16}$ & Fail & 990 & $5.496\times10^{-16}$ \\
\bottomrule
\end{tabular}
\end{table}

\subsubsection{Elliptic difference operators}
\label{subsec:neumann}
Lastly, we consider finite-difference matrices arising from discretizing second order linear-elliptic differential operators. In free space, these operators are closely related to the positive-definite kernels considered in the previous sections. We consider, for $m \in \mathbb{N}$ and $\nu + d/2\in \mathbb{N}$ and $\kappa \in \mathbb{R}$:
\begin{align}
    \begin{gathered}
        \mathcal{L}_1 := (-\Delta)^m, \quad \mathcal{L}_2 := (\kappa - \Delta)^{\nu + d/2}. 
    \end{gathered}
\end{align}
They are referred to as the polyharmonic and Matern operators, respectively. In addition to mathematical physics \cite{gilbarg1983elliptic}, these operators play an important role in statistical estimation using Gaussian random fields \cite{lindgren2011explicit}. A recurring computational kernel in these applications is the estimation of the Cholesky factors of their corresponding difference matrices. 

However, discretizing these free-space PDEs requires us to make some assumptions about the boundary conditions of these operators \cite{khristenko2019analysis}. Certain choices such as periodic or Neumann boundaries in certain domains can render the finite-difference matrix singular and in the case of certain higher-order operators (larger $m$ and $\nu$) numerically indefinite. We will demonstrate our method in this regime. 

A second-order three-point stencil in one dimension is used to discretize the Laplacian in $\mathcal{L}_1, \mathcal{L}_2$ subject to periodic and Neumann boundary conditions. The relative reconstruction accuracies of matrices corresponding to different $m$ (correspondingly $\nu + d/2$) using the finite-time methods and the pivoted Cholesky methods are compared in \cref{tab:elliptic-reconstruction-unblocked} and \cref{tab:elliptic-reconstruction-blocked}. 

\begin{table}[htbp]
\centering
\caption{Relative Frobenius reconstruction errors for 1D elliptic difference operators, \(n=500\), \(h=1/n\), normalized individually by their infinity norms. The unblocked FT-Cholesky method uses \(T=10^{14}\). Polyharmonic operators are \(L_h^m\); Mat\'ern operators are \((10^2I+L_h)^p\), where \(p=\nu+1/2\). The Neumann realization uses symmetric zero-flux differences. For Mat\'ern cases, subtract \(\delta I\) after normalization, with \(\delta=10^{2p}/s+10^{-14}\), \(s=\|A\|_\infty\), creating slightly indefinite stress tests. Comparisons use \texttt{DPOTF2} and \texttt{DPSTF2}, with default pivot tolerance.}
\label{tab:elliptic-reconstruction-unblocked}
\begin{tabular}{llrrlrr}
\toprule
Operator & Boundary & Power & FT error & \texttt{DPOTF2} & Pivoted rank & Pivoted error \\
\midrule
Polyharmonic & Neumann & 1 & $3.653\times10^{-16}$ & Fail & 499 & $1.526\times10^{-15}$ \\
 &  & 2 & $4.277\times10^{-16}$ & Pass & 499 & $2.941\times10^{-16}$ \\
 &  & 3 & $4.973\times10^{-16}$ & Fail & 500 & $1.744\times10^{-16}$ \\
\midrule
Polyharmonic & Periodic & 1 & $3.793\times10^{-16}$ & Fail & 499 & $1.559\times10^{-15}$ \\
 &  & 2 & $4.393\times10^{-16}$ & Pass & 499 & $4.272\times10^{-16}$ \\
 &  & 3 & $4.524\times10^{-16}$ & Fail & 500 & $1.751\times10^{-16}$ \\
\midrule
Mat\'ern & Neumann & 1 & $4.056\times10^{-15}$ & Fail & 499 & $1.252\times10^{-15}$ \\
 &  & 2 & $4.053\times10^{-15}$ & Fail & 499 & $2.785\times10^{-15}$ \\
 &  & 3 & $5.756\times10^{-15}$ & Fail & 499 & $6.009\times10^{-15}$ \\
\midrule
Mat\'ern & Periodic & 1 & $3.978\times10^{-15}$ & Fail & 499 & $1.260\times10^{-15}$ \\
 &  & 2 & $4.188\times10^{-15}$ & Fail & 499 & $2.539\times10^{-15}$ \\
 &  & 3 & $5.682\times10^{-15}$ & Fail & 499 & $5.726\times10^{-15}$ \\
\bottomrule
\end{tabular}
\end{table}

\begin{table}[htbp]
\centering
\caption{Relative Frobenius reconstruction errors for 1D elliptic difference operators, \(n=500\), \(h=1/n\), normalized individually by their infinity norms. The blocked FT-Cholesky method uses \(T=10^{14}\) and block size \(b=64\). Polyharmonic operators are \(L_h^m\); Mat\'ern operators are \((10^2I+L_h)^p\), where \(p=\nu+1/2\). The Neumann realization uses symmetric zero-flux differences. For Mat\'ern cases, subtract \(\delta I\) after normalization, with \(\delta=10^{2p}/s+10^{-16}\), \(s=\|A\|_\infty\), creating slightly indefinite stress tests. Comparisons use \texttt{DPOTRF} and \texttt{DPSTRF}, with default pivot tolerance.}
\label{tab:elliptic-reconstruction-blocked}
\begin{tabular}{llrrlrr}
\toprule
Operator & Boundary & Power & FT error & \texttt{DPOTRF} & Pivoted rank & Pivoted error \\
\midrule
Polyharmonic & Neumann & 1 & $3.653\times10^{-16}$ & Fail & 499 & $1.508\times10^{-15}$ \\
 &  & 2 & $4.486\times10^{-16}$ & Fail & 499 & $6.552\times10^{-16}$ \\
 &  & 3 & $5.549\times10^{-16}$ & Pass & 499 & $1.361\times10^{-15}$ \\
\midrule
Polyharmonic & Periodic & 1 & $3.793\times10^{-16}$ & Fail & 499 & $1.514\times10^{-15}$ \\
 &  & 2 & $4.111\times10^{-16}$ & Fail & 499 & $7.078\times10^{-16}$ \\
 &  & 3 & $4.150\times10^{-16}$ & Pass & 499 & $1.510\times10^{-15}$ \\
\midrule
Mat\'ern & Neumann & 1 & $4.056\times10^{-15}$ & Fail & 499 & $1.802\times10^{-15}$ \\
 &  & 2 & $4.000\times10^{-15}$ & Fail & 499 & $2.843\times10^{-15}$ \\
 &  & 3 & $5.662\times10^{-15}$ & Fail & 499 & $6.049\times10^{-15}$ \\
\midrule
Mat\'ern & Periodic & 1 & $3.978\times10^{-15}$ & Fail & 499 & $1.805\times10^{-15}$ \\
 &  & 2 & $4.183\times10^{-15}$ & Fail & 499 & $2.869\times10^{-15}$ \\
 &  & 3 & $5.886\times10^{-15}$ & Fail & 499 & $6.069\times10^{-15}$ \\
\bottomrule
\end{tabular}
\end{table}

Taken together, these three sets of experiments strongly suggest that our finite-time method is a viable alternative to pivoted Cholesky implementations. A comprehensive HPC benchmark in the spirit of \cite{hammarling2007lapack} to compare the latencies incurred due to pivoting and the associated gains made by the finite-time methods on modern computing architectures is an interesting direction for future work. 

\section{Discussion}
\label{sec:discussion}

This paper introduced the idea of finite-time regularization for factorizing positive semi-definite matrices. Canonical Cholesky factorization maps are not well defined in this regime. Therefore, the standard practice is to permute such matrices before applying this map. However, permutation incurs additional overhead, particularly for large matrices. Therefore, it is desirable to develop methods that don't resort to pivoting. To this end, we propose finite-time Cholesky factorization. 

The central idea is to design scalar dynamical systems whose finite-time evolution until $T$ provides a modified Cholesky map that is well-defined in the semi-definite case. The ambiguity of how to choose $T$ is clarified by deriving convergence rates for $\|A - L_T L_T^\top\|_F$ with respect to the evolution time $T$. 
Additionally, we have demonstrated with several numerical examples that finite-time methods are competitive with pivoted Cholesky methods. This provides an initial demonstration that the dynamical system interpretation provides a more stable Cholesky factorization method. Naturally, several question remain open. 

Firstly, to facilitate a concise presentation, we focused on the properties of the proposed Cholesky map in infinite precision. Furthermore, all numerical experiments were carried out in
dense double-precision floating-point representations of matrices. This leaves behind the obvious problem of establishing backwards stability of the finite-time Cholesky algorithm for a given floating-point representation. Naturally, one can also study how this algorithm behaves in mixed precision \cite{Higham_Mary_2022}. Understanding these questions will help us design more robust Cholesky methods in the future. 

From a practical standpoint, when the matrices are sparse (as they often are), it is not well understood whether the finite-time Cholesky map will lower the number of pivots necessary to adequately factor that matrix. 
Furthermore, in cases where matrices have a low-rank structure, the finite-time Cholesky map can improve the robustness of Nystr\"om methods \cite{gittens2013revisitingnystrommethodimproved} like RP-Cholesky \cite{chen2025randomly}, particularly when matrix representations could numerically lose positive-semi-definiteness. Empirical and theoretical investigations surrounding this can be consequential for several related problems in machine learning \cite{jalalian2025data, diaz2026robust, chen2025sparse}. 

Lastly, establishing how the finite-time error and the numerical round-off errors induced during the computation of the Cholesky map; affect downstream applications \cite{doi:10.1137/1.9780898718027.ch10} such as those where we solve linear equations, estimate the determinant of matrices, and similar matrix functions remains an active direction for future work.

\section*{Acknowledgments}
RM was partially supported by the Wallenberg AI, Autonomous Systems and
Software Program (WASP) funded by the Knut and Alice Wallenberg Foundation.

\section*{Use of AI}
The correctness of grammar in the text was verified with Grammarly.
The numerical experiments were carried out initially by RM and ChatGPT Sol + Codex was used to wrap the initial research code into a working Julia library. 
The authors verified the correctness of the code in this library. 
The proofs for Theorems 2.1, 2.2 were originally derived by RM; ChatGPT Sol was used to simplify expressions, derive limits and polish the structure of the proof. 

\appendix
\bibliographystyle{siamplain}
\bibliography{references}
\end{document}